\catcode`@=11
\expandafter\ifx\csname capitolo\endcsname\relax
      {\message {lettura formati e definizioni
standard}}\else{\endinput}\fi
\magnification \magstep1
\tolerance=1500\frenchspacing
%
%

\font\titfnt=cmbx10 scaled \magstep1
\font\pagfnt=cmsl10
\font\parfnt=cmbx10 scaled \magstep1
\font\sprfnt=cmsl10 scaled \magstep1
\font\teofnt=cmssbx10
\font\tenenufnt=cmsl10
\font\nineenufnt=cmsl9
\font\eightenufnt=cmsl8

\font\ninerm=cmr9
\font\eightrm=cmr8
\font\sixrm=cmr6
 
\font\ninei=cmmi9
\font\eighti=cmmi8
\font\sixi=cmmi6
\skewchar\ninei='177 \skewchar\eighti='177
\skewchar\sixi='177
 
\font\ninesy=cmsy9
\font\eightsy=cmsy8
\font\sixsy=cmsy6
\skewchar\ninesy='60 \skewchar\eightsy='60
\skewchar\sixsy='60

\font\ninebf=cmbx9
\font\eightbf=cmbx8
\font\sixbf=cmbx6
\font\fivebf=cmbx5
 
\font\ninett=cmtt9
\font\eighttt=cmtt8
 
\hyphenchar\tentt=-1 
\hyphenchar\ninett=-1
\hyphenchar\eighttt=-1
 
\font\ninesl=cmsl9
\font\eightsl=cmsl8
 
\font\nineit=cmti9
\font\eightit=cmti8
 
\font\tenmat=cmss10

\font\sevenmat=cmss7 
 
\font\fivemat=cmss5
\newfam\matfam

\newskip\ttglue
%
%
\newdimen\pagewidth \newdimen\pageheight
\newdimen\ruleht
\hsize=31.6pc  \vsize=44.5pc  \maxdepth=2.2pt
\parindent=19pt
\parindent=19pt
\hfuzz=2pt
\pagewidth=\hsize \pageheight=\vsize \ruleht=.5pt
\abovedisplayskip=6pt plus 3pt minus 1pt
\belowdisplayskip=6pt plus 3pt minus 1pt
\abovedisplayshortskip=0pt plus 3pt
\belowdisplayshortskip=4pt plus 3pt
\def\blank{\vskip 12pt}
\def\blankii{\blank\blank}

\def\blankq{\vskip 3pt}
%
%
\def\tenpoint{\def\rm{\fam0\tenrm}%
  \def\enufnt{\tenenufnt}%
  \textfont0=\tenrm \scriptfont0=\sevenrm
\scriptscriptfont0=\fiverm
  \textfont1=\teni \scriptfont1=\seveni
\scriptscriptfont1=\fivei
  \textfont2=\tensy \scriptfont2=\sevensy
\scriptscriptfont2=\fivesy
  \textfont3=\tenex \scriptfont3=\tenex
\scriptscriptfont3=\tenex
  \def\it{\fam\itfam\tenit}%
  \textfont\itfam=\tenit
  \def\sl{\fam\slfam\tensl}%
  \textfont\slfam=\tensl
  \def\bf{\fam\bffam\tenbf}%
  \textfont\bffam=\tenbf \scriptfont\bffam=\sevenbf
   \scriptscriptfont\bffam=\fivebf
 \def\mat{\fam\matfam\tenmat} \textfont\matfam=\tenmat
     \scriptfont\matfam=\sevenmat \scriptscriptfont\matfam=\fivemat
  \def\tt{\fam\ttfam\tentt}%
  \textfont\ttfam=\tentt
  \tt \ttglue=.5em plus.25em minus.15em
  \normalbaselineskip=12pt
  \let\sc=\eightrm
  \let\big=\tenbig
  \setbox\strutbox=\hbox{\vrule height8.5pt depth3.5pt width\z@}%
  \normalbaselines\rm}
 
\def\ninepoint{\def\rm{\fam0\ninerm}%
  \def\enufnt{\nineenufnt}%
  \textfont0=\ninerm \scriptfont0=\sixrm
\scriptscriptfont0=\fiverm
  \textfont1=\ninei \scriptfont1=\sixi
\scriptscriptfont1=\fivei
  \textfont2=\ninesy \scriptfont2=\sixsy
\scriptscriptfont2=\fivesy
  \textfont3=\tenex \scriptfont3=\tenex
\scriptscriptfont3=\tenex
  \def\it{\fam\itfam\nineit}%
  \textfont\itfam=\nineit
  \def\sl{\fam\slfam\ninesl}%
  \textfont\slfam=\ninesl
  \def\bf{\fam\bffam\ninebf}%
  \textfont\bffam=\ninebf \scriptfont\bffam=\sixbf
   \scriptscriptfont\bffam=\fivebf
  \def\tt{\fam\ttfam\ninett}%
  \textfont\ttfam=\ninett
  \tt \ttglue=.5em plus.25em minus.15em
  \normalbaselineskip=11pt
  \let\sc=\sevenrm
  \let\big=\ninebig
  \setbox\strutbox=\hbox{\vrule height8pt depth3pt width\z@}%
  \normalbaselines\rm}
 
\def\eightpoint{\def\rm{\fam0\eightrm}%
  \def\enufnt{\eightenufnt}%
  \textfont0=\eightrm \scriptfont0=\sixrm
\scriptscriptfont0=\fiverm
  \textfont1=\eighti \scriptfont1=\sixi
\scriptscriptfont1=\fivei
  \textfont2=\eightsy \scriptfont2=\sixsy
\scriptscriptfont2=\fivesy
  \textfont3=\tenex \scriptfont3=\tenex
\scriptscriptfont3=\tenex
  \def\it{\fam\itfam\eightit}%
  \textfont\itfam=\eightit
  \def\sl{\fam\slfam\eightsl}%
  \textfont\slfam=\eightsl
  \def\bf{\fam\bffam\eightbf}%
  \textfont\bffam=\eightbf \scriptfont\bffam=\sixbf
   \scriptscriptfont\bffam=\fivebf
  \def\tt{\fam\ttfam\eighttt}%
  \textfont\ttfam=\eighttt
  \tt \ttglue=.5em plus.25em minus.15em
  \normalbaselineskip=9pt
  \let\sc=\sixrm
  \let\big=\eightbig
  \setbox\strutbox=\hbox{\vrule height7pt depth2pt width\z@}%
  \normalbaselines\rm}
 
\def\tenbig#1{{\hbox{$\left#1\vbox to8.5pt{}\right.\n@space$}}}
\def\ninebig#1{{\hbox{$\textfont0=\tenrm\textfont2=\tensy
  \left#1\vbox to7.25pt{}\right.\n@space$}}}
\def\eightbig#1{{\hbox{$\textfont0=\ninerm\textfont2=\ninesy
  \left#1\vbox to6.5pt{}\right.\n@space$}}}
%
%
\def\corsivo{\enufnt}

\def\footnote#1{\edef\@sf{\spacefactor\the\spacefactor}$^{#1}$\@sf
      \insert\footins\bgroup\ninepoint
      \interlinepenalty100 \let\par=\endgraf
        \leftskip=\z@skip \rightskip=\z@skip
        \splittopskip=10pt plus 1pt minus 1pt \floatingpenalty=20000
        \smallskip\item{$^{#1}$}\bgroup\strut\aftergroup\@foot\let\next}
\skip\footins=12pt plus 2pt minus 4pt 
\dimen\footins=30pc 
%
%
\newcount\numbibliogr@fi@
\global\numbibliogr@fi@=1
\newwrite\fileref
\immediate\openout\fileref=\jobname.ref
\immediate\write\fileref{\parindent 30pt}
\def\cita#1#2{\def\us@gett@{\the\numbibliogr@fi@}
\expandafter\xdef\csname:bib_#1\endcsname{\us@gett@}
           \immediate\write\fileausiliario{\string\expandafter\string\edef
           \string\csname:bib_#1\string\endcsname{\us@gett@}}
           \immediate\write\fileref
           {\par\noexpand\item{{[\the\numbibliogr@fi@]\enspace}}}\ignorespaces
           \immediate\write\fileref{{#2}}\ignorespaces
           \global\advance\numbibliogr@fi@ by 1\ignorespaces}
\def\bibref#1{\seindefinito{:bib_#1}
          \immediate\write16{ !!! \string\bibref{#1} non definita !!!}
 
\expandafter\xdef\csname:bib_#1\endcsname{??}\fi
      {$^{[\csname:bib_#1\endcsname]}$}}
\def\dbiref#1{\seindefinito{:bib_#1}
          \immediate\write16{ !!! \string\bibref{#1} non definita !!!}
 
\expandafter\xdef\csname:bib_#1\endcsname{??}\fi
      {[\csname:bib_#1\endcsname]}}

\def\references{\immediate\closeout\fileref
                \par\goodbreak
                \blankii
                \centerline{\parfnt References}
                \nobreak\blank\nobreak
                \input \jobname.ref}

%
%
\def\title#1{\null\blankii\noindent{\titfnt\uppercase{#1}}\blank}
\def\riga{\par\vskip 6pt\noindent}
\def\author#1{\leftskip 1.8cm\smallskip\noindent{#1}\smallskip\leftskip 0pt}
\def\abstract#1{\par\blankii\noindent
          {\ninepoint
            {\bf Abstract. }{\rm #1}}
          \par}

\def\sunto#1{\par\blankii\noindent
          {\ninepoint
            {\bf Sunto. }{\rm #1}}
          \par}
%
%
%
\def\seindefinito#1{\expandafter\ifx\csname#1\endcsname\relax}
\newwrite\filesimboli
    \immediate\openout\filesimboli=\jobname.smb
\newwrite\fileausiliario
    \openin1 \jobname.aux
 
\ifeof1\relax\else\closein1\relax\input\jobname.aux\fi
    \immediate\openout\fileausiliario=\jobname.aux
%
\newdimen\@mpiezz@
\@mpiezz@=\hsize
\newbox\boxp@r@gr@fo
\newcount\nump@r@gr@fo
\global\nump@r@gr@fo=0
\def\section#1#2{
           \global\advance\nump@r@gr@fo by 1
           \xdef\paragrafocorrente{\the\nump@r@gr@fo}
\expandafter\xdef\csname:sec_#1\endcsname{\paragrafocorrente}
           \global\numsottop@r@gr@fo=0
      \immediate\write\filesimboli{  }
      \immediate\write\filesimboli{  }
      \immediate\write\filesimboli{--> Section \paragrafocorrente :
               rif.: #1}
\immediate\write\fileausiliario{\string\expandafter\string\edef
\string\csname:sec_#1\string\endcsname{\paragrafocorrente}}
      \write16{Section {\csname:sec_#1\endcsname}.}
           \goodbreak\vskip 24pt plus 6pt\noindent\ignorespaces
              {\setbox\boxp@r@gr@fo=\hbox{\parfnt\noindent\ignorespaces
              {\paragrafocorrente.\quad}}\ignorespaces
             \advance\@mpiezz@ by -\wd\boxp@r@gr@fo
             \box\boxp@r@gr@fo\vtop{\hsize=\@mpiezz@\noindent\parfnt #2}}
           \par\nobreak\vskip 12pt plus 3pt\nobreak
           \noindent\ignorespaces}
\def\secref#1{\seindefinito{:sec_#1}
          \immediate\write16{ !!! \string\secref{#1} non definita !!!}
 \expandafter\xdef\csname:sec_#1\endcsname{??}\fi
      \csname:sec_#1\endcsname
      }
\newcount\num@ppendice
\global\num@ppendice=64
\def\appendix#1#2{
           \global\advance\num@ppendice by 1
           \xdef\paragrafocorrente{\char\the\num@ppendice}
           \expandafter\xdef\csname:app_#1\endcsname{\paragrafocorrente}
           \global\numsottop@r@gr@fo=0
      \immediate\write\filesimboli{  }
      \immediate\write\filesimboli{  }
      \immediate\write\filesimboli{--> Appendice \paragrafocorrente :
               rif.: #1}
      \immediate\write\fileausiliario{\string\expandafter\string\edef
      \string\csname:app_#1\string\endcsname{\paragrafocorrente}}
      \write16{Appendice {\csname:app_#1\endcsname}.}
           \goodbreak\vskip 18pt plus 6pt\noindent\ignorespaces
              {\setbox\boxp@r@gr@fo=\hbox{\parfnt\noindent\ignorespaces
              {\paragrafocorrente.\quad}}\ignorespaces
             \advance\@mpiezz@ by -\wd\boxp@r@gr@fo
             \box\boxp@r@gr@fo\vtop{\hsize=\@mpiezz@\noindent\parfnt #2}}
           \par\nobreak\vskip 9pt plus 3pt\nobreak
           \noindent\ignorespaces}
\def\appref#1{\seindefinito{:app_#1}
          \immediate\write16{ !!! \string\appref{#1} non definita !!!}
          \expandafter\xdef\csname:app_#1\endcsname{??}\fi
      \csname:app_#1\endcsname
      }
%
\newcount\numsottop@r@gr@fo
\global\numsottop@r@gr@fo=0
\def\subsection#1#2{
           \global\advance\numsottop@r@gr@fo by 1
 \xdef\sottoparcorrente{\paragrafocorrente.\the\numsottop@r@gr@fo}
 \expandafter\xdef\csname:sbs_#1\endcsname{\sottoparcorrente}
      \immediate\write\filesimboli{  }
      \immediate\write\filesimboli{-----> Subsection \sottoparcorrente :
               rif.: #1}
 \immediate\write\fileausiliario{\string\expandafter\string\edef
 \string\csname:sbs_#1\string\endcsname{\sottoparcorrente}}
      \write16{Subsection {\csname:sbs_#1\endcsname}.}
           \goodbreak\vskip 9pt plus 2pt\noindent\ignorespaces
 {\setbox\boxp@r@gr@fo=\hbox{\sprfnt\noindent\ignorespaces
              {\sottoparcorrente\quad}}\ignorespaces
             \advance\@mpiezz@ by -\wd\boxp@r@gr@fo
 \box\boxp@r@gr@fo\vtop{\hsize=\@mpiezz@\noindent\sprfnt#2}}
           \par\nobreak\vskip 3pt plus 1pt\nobreak
           \noindent\ignorespaces}
\def\sbsref#1{\seindefinito{:sbs_#1}
          \immediate\write16{ !!! \string\sbsref{#1} non definita !!!}
 \expandafter\xdef\csname:sbs_#1\endcsname{??}\fi
      \csname:sbs_#1\endcsname
      }
%

\newcount\numformul@
\global\numformul@=0
\def\formula#1{
            \formdef{#1}\ignorespaces
            \leqno{(\csname:frm_#1\endcsname)}
            }
\def\formdef#1{
            \global\advance\numformul@ by 1
            \def\us@gett@{\the\numformul@}
 \expandafter\xdef\csname:frm_#1\endcsname{\us@gett@}
      \immediate\write\filesimboli{--------> Formula \us@gett@ :
               rif.: #1}
 \immediate\write\fileausiliario{\string\expandafter\string\edef
 \string\csname:frm_#1\string\endcsname{\us@gett@}}
            }
\def\frmref#1{\seindefinito{:frm_#1}
          \immediate\write16{ !!! \string\frmref{#1} non definita !!!}
 \expandafter\xdef\csname:frm_#1\endcsname{??}\fi
      (\csname:frm_#1\endcsname)}
%
\newcount\numenunci@to
\global\numenunci@to=0
\def\endclaim{\endgroup
            \par\if F\sp@zi@tur@{\blankq}\gdef\sp@zi@tur@{T}\fi}
\newcount\numteorema
\global\numteorema=0
\def\theorem#1{
            \global\advance\numteorema by 1
            \def\us@gett@{\the\numteorema}
 \expandafter\xdef\csname:thr_#1\endcsname{\us@gett@}
      \immediate\write\filesimboli{--------> Theorem \us@gett@ :
               rif.: #1}
 \immediate\write\fileausiliario{\string\expandafter\string\edef
 \string\csname:thr_#1\string\endcsname{\us@gett@}}
            \par\if T\sp@zi@tur@{\gdef\sp@zi@tur@{F}}\else{\blankq}\fi
            \noindent{\teofnt Theorem \us@gett@:\quad}\begingroup\enufnt
            \ignorespaces}
\def\theoremnn{
            \par\if T\sp@zi@tur@{\gdef\sp@zi@tur@{F}}\else{\blankq}\fi
            \noindent{\teofnt Theorem:\quad}\begingroup\enufnt
            \ignorespaces}
\def\theoremtx#1#2{
            \global\advance\numteorema by 1
            \def\us@gett@{\the\numteorema}
 \expandafter\xdef\csname:thr_#1\endcsname{\us@gett@}
      \immediate\write\filesimboli{--------> Theorem \us@gett@ :
               rif.: #1}
 \immediate\write\fileausiliario{\string\expandafter\string\edef
 \string\csname:thr_#1\string\endcsname{\us@gett@}}
            \par\if T\sp@zi@tur@{\gdef\sp@zi@tur@{F}}\else{\blankq}\fi
            \noindent{\teofnt Theorem \us@gett@:
                {\enufnt #2}.\quad}\begingroup\enufnt
            \ignorespaces}
\def\thrref#1{\seindefinito{:thr_#1}
          \immediate\write16{ !!! \string\thrref{#1} non definita !!!}
 \expandafter\xdef\csname:thr_#1\endcsname{??}\fi
      \csname:thr_#1\endcsname}
%
\newcount\numproposizione
\global\numproposizione=0
\def\proposition#1{
            \global\advance\numproposizione by 1
            \def\us@gett@{\the\numproposizione}
 \expandafter\xdef\csname:pro_#1\endcsname{\us@gett@}
      \immediate\write\filesimboli{--------> Proposition \us@gett@ :
               rif.: #1}
 \immediate\write\fileausiliario{\string\expandafter\string\edef
 \string\csname:pro_#1\string\endcsname{\us@gett@}}
            \par\if T\sp@zi@tur@{\gdef\sp@zi@tur@{F}}\else{\blankq}\fi
            \noindent{\teofnt Proposition \us@gett@:\quad}\begingroup\enufnt
            \ignorespaces}
\def\propositionn{
            \par\if T\sp@zi@tur@{\gdef\sp@zi@tur@{F}}\else{\blankq}\fi
            \noindent{\teofnt Proposition:\quad}\begingroup\enufnt
            \ignorespaces}
\def\propositiontx#1#2{
            \global\advance\numproposizione by 1
            \def\us@gett@{\the\numproposizione}
 \expandafter\xdef\csname:pro_#1\endcsname{\us@gett@}
      \immediate\write\filesimboli{--------> Proposition \us@gett@ :
               rif.: #1}
  \immediate\write\fileausiliario{\string\expandafter\string\edef
 \string\csname:pro_#1\string\endcsname{\us@gett@}}
            \par\if T\sp@zi@tur@{\gdef\sp@zi@tur@{F}}\else{\blankq}\fi
            \noindent{\teofnt Proposition \us@gett@:
               {\enufnt #2}.\quad}\begingroup\enufnt
            \ignorespaces}
\def\proref#1{\seindefinito{:pro_#1}
          \immediate\write16{ !!! \string\proref{#1} non definita !!!}
 \expandafter\xdef\csname:pro_#1\endcsname{??}\fi
      \csname:pro_#1\endcsname}
%
\newcount\numcorollario
\global\numcorollario=0
\def\corollary#1{
            \global\advance\numcorollario by 1
            \def\us@gett@{\the\numcorollario}
 \expandafter\xdef\csname:cor_#1\endcsname{\us@gett@}
      \immediate\write\filesimboli{--------> Corollary \us@gett@ :
               rif.: #1}
 \immediate\write\fileausiliario{\string\expandafter\string\edef
 \string\csname:cor_#1\string\endcsname{\us@gett@}}
            \par\if T\sp@zi@tur@{\gdef\sp@zi@tur@{F}}\else{\blankq}\fi
            \noindent{\teofnt Corollary
\us@gett@:\quad}\begingroup\enufnt
            \ignorespaces}
\def\corref#1{\seindefinito{:cor_#1}
          \immediate\write16{ !!! \string\corref{#1} non definita !!!}
 \expandafter\xdef\csname:cor_#1\endcsname{??}\fi
      \csname:cor_#1\endcsname}
%
\newcount\numlemma
\global\numlemma=0
\def\lemma#1{
            \global\advance\numlemma by 1
            \def\us@gett@{\the\numlemma}
\expandafter\xdef\csname:lem_#1\endcsname{\us@gett@}
      \immediate\write\filesimboli{--------> Lemma \us@gett@ :
               rif.: #1}
 \immediate\write\fileausiliario{\string\expandafter\string\edef
 \string\csname:lem_#1\string\endcsname{\us@gett@}}
            \par\if T\sp@zi@tur@{\gdef\sp@zi@tur@{F}}\else{\blankq}\fi
            \noindent{\teofnt Lemma \us@gett@:\quad}\begingroup\enufnt
            \ignorespaces}
\def\lemref#1{\seindefinito{:lem_#1}
          \immediate\write16{ !!! \string\lemref{#1} non definita !!!}
 \expandafter\xdef\csname:lem_#1\endcsname{??}\fi
      \csname:lem_#1\endcsname}
%
\newcount\numdefinizione
\global\numdefinizione=0
\def\definition#1{
            \global\advance\numdefinizione by 1
            \def\us@gett@{\the\numdefinizione}
 \expandafter\xdef\csname:def_#1\endcsname{\us@gett@}
      \immediate\write\filesimboli{--------> Definition \us@gett@ :
               rif.: #1}
 \immediate\write\fileausiliario{\string\expandafter\string\edef
 \string\csname:def_#1\string\endcsname{\us@gett@}}
            \par\if T\sp@zi@tur@{\gdef\sp@zi@tur@{F}}\else{\blankq}\fi
            \noindent{\teofnt Definition \us@gett@:\quad}\begingroup\enufnt
            \ignorespaces}
\def\definitiontx#1#2{
            \global\advance\numdefinizione by 1
            \def\us@gett@{\the\numdefinizione}
 \expandafter\xdef\csname:def_#1\endcsname{\us@gett@}
      \immediate\write\filesimboli{--------> Definition \us@gett@ :
               rif.: #1}
 \immediate\write\fileausiliario{\string\expandafter\string\edef
 \string\csname:def_#1\string\endcsname{\us@gett@}}
            \par\if T\sp@zi@tur@{\gdef\sp@zi@tur@{F}}\else{\blankq}\fi
            \noindent{\teofnt Definition \us@gett@:
               {\enufnt #2}.\quad}\begingroup\enufnt
            \ignorespaces}
\def\defref#1{\seindefinito{:def_#1}
          \immediate\write16{ !!! \string\defref{#1} non definita !!!}
 \expandafter\xdef\csname:def_#1\endcsname{??}\fi
      \csname:def_#1\endcsname}
%
\def\proof{\par\if T\sp@zi@tur@{\gdef\sp@zi@tur@{F}}\else{\blankq}\fi
    \noindent{\teofnt Proof.\quad}\begingroup\ignorespaces}
\def\prooftx#1{\par\if T\sp@zi@tur@{\gdef\sp@zi@tur@{F}}\else{\blankq}\fi
    \noindent{\teofnt Proof #1.\quad}\begingroup\ignorespaces}
\def\endproof{\nobreak\quad\nobreak\hfill\nobreak{\enufnt Q.E.D.}
    \endclaim}
%
\newcount\numfigur@
\global\numfigur@=0
\newbox\boxfigur@
\newbox\comfigur@
\newdimen\@mpfigur@
\newdimen\m@rfigur@
\m@rfigur@=2 pc
\@mpfigur@=\hsize
\advance\@mpfigur@ by -2\m@rfigur@
\def\figure#1#2#3{
            \global\advance\numfigur@ by 1
            \def\us@gett@{\the\numfigur@}
            \expandafter\xdef\csname:fig_#1\endcsname{\us@gett@}
      \immediate\write\filesimboli{--------> Figure \us@gett@ :
               rif.: #1}
      \immediate\write\fileausiliario{\string\expandafter\string\edef
      \string\csname:fig_#1\string\endcsname{\us@gett@}}
      \setbox\boxfigur@\vbox{\centerline{#2}}
      \topinsert
         {\vbox{
               \vskip 1pt
               \box\boxfigur@
               \vskip 4pt}}
      \setbox\comfigur@\vtop{\hsize=\@mpfigur@\parindent 0pt
         {\ninepoint
         {\teofnt Figure \us@gett@.}\enspace{#3}}
      }
      \centerline{\box\comfigur@}
      \endinsert
      \write16{Figure {\csname:fig_#1\endcsname}.}
      }
\def\figcont#1#2{
      \setbox\boxfigur@\vbox{\centerline{#2}}
      \topinsert
         {\vbox{
               \vskip 1pt
               \box\boxfigur@
               \vskip 4pt}}
 \setbox\comfigur@\vtop{\hsize=\@mpfigur@\parindent 0pt
         {\ninepoint
         {\teofnt Figure \figref{#1}.}\enspace{(continued).}}
      }
      \centerline{\box\comfigur@}
      \endinsert
      \write16{Figura (cont) {\csname:fig_#1\endcsname}.}
}
\def\figref#1{\seindefinito{:fig_#1}
          \immediate\write16{ !!! \string\figref{#1} non definita !!!}
 \expandafter\xdef\csname:fig_#1\endcsname{??}\fi
      \csname:fig_#1\endcsname}
\def\citazione{\begingroup
     \everypar{\parshape 1 \m@rfigur@ \@mpfigur@\noindent}}

%
\def\remarks{\par\if
T\sp@zi@tur@{\gdef\sp@zi@tur@{F}}\else{\blankq}\fi
    \noindent{\teofnt Remarks.\quad}\ignorespaces}
\def\remark{\par\if
T\sp@zi@tur@{\gdef\sp@zi@tur@{F}}\else{\blankq}\fi
    \noindent{\teofnt Remark.\quad}\ignorespaces}
%
\newcount\numt@vol@
\global\numt@vol@=0
\newbox\boxt@vol@
\newbox\comt@vol@
\newdimen\@mpt@vol@
\newdimen\m@rt@vol@
\m@rt@vol@=2 pc
\@mpt@vol@=\hsize
\advance\@mpt@vol@ by -2\m@rt@vol@
\def\table#1#2#3{
            \global\advance\numt@vol@ by 1
            \def\us@gett@{\the\numt@vol@}
            \expandafter\xdef\csname:tav_#1\endcsname{\us@gett@}
      \immediate\write\filesimboli{--------> Table \us@gett@ :
               rif.: #1}
      \immediate\write\fileausiliario{\string\expandafter\string\edef
      \string\csname:tav_#1\string\endcsname{\us@gett@}}
      \setbox\boxt@vol@\vbox{#3}
      \topinsert
      \setbox\comt@vol@\vtop{\hsize=\@mpt@vol@\parindent 0pt
         {\ninepoint
         {\teofnt Table \us@gett@.}\enspace{#2}}
      }
      \centerline{\box\comt@vol@}
      \vskip 5pt%
      \centerline{\box\boxt@vol@}%
      \endinsert
      \write16{Table {\csname:tav_#1\endcsname}.}
}
\def\tabcont#1#2{
      \setbox\boxt@vol@\vbox{#2}
      \topinsert
      \setbox\comt@vol@\vtop{\hsize=\@mpt@vol@\parindent 0pt
         {\ninepoint
         {\teofnt Table \tabref{#1}.}\enspace{(continued)}}
      }
      \centerline{\box\comt@vol@}
      \vskip 5pt%
      \centerline{\box\boxt@vol@}%
      \endinsert
      \write16{Table {\csname:tav_#1\endcsname}.}
}
\def\tabref#1{\seindefinito{:tav_#1}
          \immediate\write16{ !!! \string\tabref{#1} non definita !!!}
          \expandafter\xdef\csname:tav_#1\endcsname{??}\fi
      \csname:tav_#1\endcsname}
%
%

\def\noblank{\gdef\sp@zi@tur@{T}}
\def\incolonna#1{\displ@y \tabskip=\centering
   \halign to \displaywidth{\hfil$\@lign\displaystyle{{}##}$\tabskip=0pt
       &$\@lign\displaystyle{{}##}$\hfil\tabskip=\centering
       &\llap{$\@lign\displaystyle{{}##}$}\tabskip=0pt\crcr
       #1\crcr}}
%
%
\nopagenumbers
\def\testos{\null}
\def\testod{\null}
\headline={\if T\tpage{\gdef\tpage{F}{\hfil}}
 \else{\ifodd\pageno\rightheadline\else\leftheadline\fi}
           \fi}
 
\gdef\tpage{T}
\def\rightheadline{\hfil{\pagfnt\testod}\hfil{\pagfnt\folio}}
\def\leftheadline{{\pagfnt\folio}\hfil{\pagfnt\testos}\hfil}
\voffset=2\baselineskip
\everypar={\gdef\sp@zi@tur@{F}}
\catcode`@=12
%

\def\pmb#1{\setbox0=\hbox{#1}\ignorespaces
    \hbox{\kern-.02em\copy0\kern-\wd0\ignorespaces
    \kern.05em\copy0\kern-\wd0\ignorespaces
    \kern-.02em\raise.02em\box0 }}
\def\gt{>}

\def\rho{\varrho}
\def\theta{\vartheta}
\let\flusso=\phi       
\def\phi{\varphi}
\def\epsilon{\varepsilon}

\def\fraz#1#2{{{#1}\over{#2}}}
\def\frac#1#2{{{#1}\over{#2}}}

\def\diff{\mathord{\rm d}}

\def\der#1#2{{{{\diff}{#1}}\over{{\diff}{#2}}}}

\def\derpar#1#2{{{{\partial}{#1}}\over{{\partial}{#2}}}}

\def\reali{\mathinner{\bf R}}
\def\complessi{\mathinner{\bf C}}
\def\toro{\mathinner{\bf T}}
\def\interi{\mathinner{\bf Z}}

\def\poisson#1#2{\lbrace#1,#2\rbrace}

\tenpoint\rm


\def\lie#1{L_{#1}}
\def\liec#1{L_{#1}}

\def\Dscr{{\cal D}}

\def\Gscr{{\cal G}}

\def\Nscr{{\cal N}}
\def\Oscr{{\cal O}}
\def\Pscr{{\cal P}}

\def\Rscr{{\cal R}}

\def\evet{{\bf e}}
\def\Amat{{\mat A}}
\def\Bmat{{\mat B}}
\def\Imat{{\mat I}}
\def\Jmat{{\mat J}}
\def\Dmat{{\mat D}}
\def\Rmat{{\mat R}}
\def\Rmatinv{{{\mat R}^{-1}}}
\def\Lambdamat{{\mat\Lambda}}
\def\Omegamat{{\mat\Omega}}

\def\poisson#1#2{\lbrace #1,#2 \rbrace}

\def\diag{\mathop {\rm diag}}

\def\uno{{\mat 1}}


\cita{Birkhoff-1913}{G.D.\ Birkhoff: {\it
Proof of Poincaré's Geometric Theorem}, Transactions of the American
Mathematical Society {\bf 14}, 14--22 (1913)}

\cita{Birkhoff-1920}{G.D.\ Birkhoff: {Surface transformations and their dynamical
applications}, Acta Mathematica {\bf 43}, 1--119 (1920).}

\cita{Deprit-1969}{A. Deprit: {\it Canonical transformations
depending on a small parameter}, {Cel. Mech.} {\bf 1}, 12--30 (1969).}

\cita{Fasso-1989}{F. Fass\`o: {\it On a relation among Lie series},
Celestial Mechanics {\bf 46}, 113--118 (1989).}

\cita{Giorgilli-1978}{A. Giorgilli, L. Galgani: {\it Formal
integrals for an autonomous Hamiltonian system near an equilibrium
point}, Cel. Mech. {\bf 17} (1978), 267--280.
}
\cita{Giorgilli-1985}{A. Giorgilli, L. Galgani: {\it Rigorous
estimates for the series expansions of Hamiltonian perturbation
theory}, Cel. Mech.  {\bf 37}, 95-112 (1985).}

\cita{Giorgilli-1997.2}{A. Giorgilli, U. Locatelli: {\it Kolmogorov theorem 
and classical perturbation theory}, ZAMP {\bf 48}, 220--261 (1997).}

\cita{Giorgilli-1999}{A. Giorgilli, U. Locatelli: {\it A classical
self--contained proof of Kolmogorov's theorem on invariant tori}, in
Proceedings of the NATO ASI school ``Hamiltonian systems with three or
more degrees of freedom'', C.\ Sim\'o ed., NATO ASI series 
C:\ Math.\ Phys.\ Sci., Vol.\ 533, Kluwer Academic Publishers,
Dordrecht--Boston--London, 72--89 (1999).}

\cita{Giorgilli-2010}{A.\ Giorgilli,  S.\ Marmi: {\it Improved estimates for the
convergence radius in the Poincar\'e--Siegel problem}, Discrete and
Continuous Dynamical Systems series S 3, 601--621 (2010).}

\cita{Giorgilli-2012.1}{A. Giorgilli, M. Sansottera: {\it Methods of
algebraic manipulation in perturbation theory}, in {\it Chaos,
Diffusion and Non-integrability in Hamiltonian Systems - Applications
to Astronomy}, Proceedings of the 3rd La Plata International School on
Astronomy and Geophysics, P.M. Cincotta, C.M. Giordano and
C. Efthymiopoulos eds., Universidad Nacional de La Plata and
Asociación Argentina de Astronomía Publishers, La Plata, Argentina
(2012). }

\cita{Groebner-1957}{W.\ Gr\"obner: {\it Nuovi contributi alla 
teoria dei sistemi di equazioni differenziali nel campo analitico},
Atti Accad. Naz. Lincei. Rend. Cl. Sci. Fis. Mat. Nat. {\bf 23}
375-–379 (1957).}

\cita{Groebner-1960}{W.\ Gr\"obner: {\it Die  Lie--Reihen und Ihre
Anwendungen}, VEB Deutscher Verlag der Wissenschaften, Mathematische 
Monographien, {\bf 3} (1960).
Italian translation: {\it Serie di Lie e loro applicazioni},
Ed. Cremonese, Roma, 1973.}

\cita{Haro-2000}{A.\ Haro: {\it The primitive function of an exact
symplectomorphism}, Nonlinearity {\bf 13}, 1483--1500 (2000).}

\cita{Henon-1969}{M.\ H\'enon: {\it Numerical Study of Quadratic
Area-Preserving Mappings}, Quart. Appl. Math. {\bf 27}, 291--312
(1969).}

\cita{Henrard-1974}{J.\ Henrard, J.\ Roels: {\it Equivalence for Lie
transforms}, Celestial Mechanics {\bf 10}, 497--512 (1974).}

\cita{Hori-1966}{G. Hori: {\it Theory of general perturbations with unspecified
canonical variables}, {Publ. Astron. Soc. Japan}, {\bf 18},
287--296  (1966).}


\cita{Kuksin-1993}{S.\ Kuksin: {\it    On the Inclusion of an Analytic
Symplectomorphism Close to an Integrable One into a Hamiltonian Flow}, Russian
Journal of Mathematical Physics {\bf 1}, 191--207 (1993).}

\cita{Kuksin-1994}{S.\ Kuksin, J.\ P\"oschel: {\it On the inclusion of analytic
symplectic maps in analytic Hamiltonian flows and its applications}, in: S.\
Kuksin, V.\ Lazutkin, J\ P\"oschel (Eds.), {\it
Seminar on dynamical systems}, Birkh\"auser, Basel, 96--116 (1994)}

\cita{Moser-1962}{J.K.\ Moser: {\it On invariant curves of area--preserving 
mappings of an annulus}, Nachr.\ Akad.\ Wiss.\ G\"ott., II Math.\
Phys.\ Kl.\ 1962, 1--20 (1962).}

\cita{Poincare-1892}{Poincar\'e, H.: {\it Les m\'ethodes nouvelles de la 
m\'ecanique c\'eleste}, Gauthier--Villars, Paris (1892).}

\cita{Poincare-1912}{H.\ Poincar\'e: {\it Sur un th\'eor\`eme de
g\`eom\'etrie}, Rendiconti del Circolo Matematico di Palermo {\bf
33}, 375--407 (1912).} 

\cita{Pronin-1997}{A.V.\ Pronin, D.V.\ Treschev: {\it On the inclusion of
analytic maps into analytic flows}, Regular and Chaotic Dynamics {\bf 2}, 14--24
(1997).}

\cita{Schroeder-1871}{E. Schr\"oder: {\it \"Uber iterierte
Functionen}, Math. Ann. {\bf 3}, 296-322 (1871).}

\cita{Siegel-1942}{Siegel, C.L.: {\it Iterations of analytic functions},
Annals of Math. {\bf 43}, 607--612 (1942).}

\def\testos{A. Giorgilli}
\def\testod{On the representation of maps by Lie transforms}

\expandafter\edef\csname:sec_liemaps.1\endcsname{1}
\expandafter\edef\csname:sec_liemaps.2\endcsname{2}
\expandafter\edef\csname:sbs_liemaps.2.1\endcsname{2.1}
\expandafter\edef\csname:frm_lie.1\endcsname{1}
\expandafter\edef\csname:frm_lie.2\endcsname{2}
\expandafter\edef\csname:frm_lie.3\endcsname{3}
\expandafter\edef\csname:frm_diffeq.6\endcsname{4}
\expandafter\edef\csname:sbs_liemaps.2.2\endcsname{2.2}
\expandafter\edef\csname:frm_lser.1\endcsname{5}
\expandafter\edef\csname:frm_lser.2\endcsname{6}
\expandafter\edef\csname:frm_lser.3\endcsname{7}
\expandafter\edef\csname:frm_trslie.1\endcsname{8}
\expandafter\edef\csname:frm_trslie.2\endcsname{9}
\expandafter\edef\csname:frm_trslie.5\endcsname{10}
\expandafter\edef\csname:frm_trslie.6\endcsname{11}
\expandafter\edef\csname:frm_lser.6\endcsname{12}
\expandafter\edef\csname:frm_trslie.3\endcsname{13}
\expandafter\edef\csname:frm_trslie.4\endcsname{14}
\expandafter\edef\csname:sbs_liemaps.2.3\endcsname{2.3}
\expandafter\edef\csname:frm_lser.7\endcsname{15}
\expandafter\edef\csname:pro_trslie.10\endcsname{1}
\expandafter\edef\csname:frm_trslie.7\endcsname{16}
\expandafter\edef\csname:frm_lser.8r\endcsname{17}
\expandafter\edef\csname:frm_lser.8\endcsname{18}
\expandafter\edef\csname:frm_lser.9\endcsname{19}
\expandafter\edef\csname:pro_lser.10\endcsname{2}
\expandafter\edef\csname:sbs_liemaps.2.4\endcsname{2.4}
\expandafter\edef\csname:frm_trslie.8\endcsname{20}
\expandafter\edef\csname:pro_trslie.11\endcsname{3}
\expandafter\edef\csname:frm_trslie.12\endcsname{21}
\expandafter\edef\csname:frm_trslie.13a\endcsname{22}
\expandafter\edef\csname:sec_liemaps.3\endcsname{3}
\expandafter\edef\csname:frm_map.1\endcsname{23}
\expandafter\edef\csname:frm_map.2\endcsname{24}
\expandafter\edef\csname:sbs_liemaps.3.1\endcsname{3.1}
\expandafter\edef\csname:frm_map.1a\endcsname{25}
\expandafter\edef\csname:frm_map.3\endcsname{26}
\expandafter\edef\csname:frm_map.2a\endcsname{27}
\expandafter\edef\csname:sbs_liemaps.3.2\endcsname{3.2}
\expandafter\edef\csname:frm_map.4\endcsname{28}
\expandafter\edef\csname:pro_map.13\endcsname{4}
\expandafter\edef\csname:frm_map.5\endcsname{29}
\expandafter\edef\csname:sbs_liemaps.3.3\endcsname{3.3}
\expandafter\edef\csname:frm_map.6\endcsname{30}
\expandafter\edef\csname:frm_map.7\endcsname{31}
\expandafter\edef\csname:pro_map.10\endcsname{5}
\expandafter\edef\csname:frm_map.9\endcsname{32}
\expandafter\edef\csname:frm_map.20c\endcsname{33}
\expandafter\edef\csname:frm_map.20b\endcsname{34}
\expandafter\edef\csname:pro_map.20\endcsname{6}
\expandafter\edef\csname:frm_map.21a\endcsname{35}
\expandafter\edef\csname:frm_map.21b\endcsname{36}
\expandafter\edef\csname:frm_map.21c\endcsname{37}
\expandafter\edef\csname:sec_liemaps.4\endcsname{4}
\expandafter\edef\csname:sbs_liemaps.4.1\endcsname{4.1}
\expandafter\edef\csname:frm_nrm.1\endcsname{38}
\expandafter\edef\csname:sbs_liemaps.4.2\endcsname{4.2}
\expandafter\edef\csname:frm_nrm.2\endcsname{39}
\expandafter\edef\csname:frm_nrm.3\endcsname{40}
\expandafter\edef\csname:frm_map.22\endcsname{41}
\expandafter\edef\csname:frm_map.23\endcsname{42}
\expandafter\edef\csname:app_liemaps.a\endcsname{\char 65}
\expandafter\edef\csname:frm_trslie.5b\endcsname{43}
\expandafter\edef\csname:frm_trslie.13\endcsname{44}
\expandafter\edef\csname:frm_trslie.8a\endcsname{45}
\expandafter\edef\csname:frm_trslie.8b\endcsname{46}

\noindent
To appear in: \hfill\break
{\sl Rendiconti dell'Istituto Lombardo Accademia di Scienze e Lettere,
Classe di Scienze.}
\vskip 24 pt

\title{On the representation of maps\riga by Lie transforms}

\author{\it ANTONIO GIORGILLI
\hfill\break Dipartimento di Matematica,
Via Saldini 50,
\hfill\break 20133\ ---\  Milano, Italy.}

\sunto{Si riconsidera il problema della rappresentazione di una mappa
in una forma adatta all'applicazione dei metodi di forma normale.  Si
mostra che ricorrendo ai metodi delle serie di Lie e delle trasformate
di Lie si pu\`o costruire in modo diretto un algoritmo di
normalizzazione.  Si discute brevemente l'applicazione alla mappa di
Schr\"oder--Siegel e alla mappa standard di Chirikov, estendondole al
caso di dimensione generica.}

\abstract{The problem of representing a class of maps in a form suited for
application of normal form methods is revisited.  It is shown that
using the methods of Lie series and of Lie transform a normal form
algorithm is constructed in a straightforward manner.  The examples of
the Schr\"oder--Siegel map and of the Chirikov standard map are included, with
extension to arbitrary dimension.}

\section{liemaps.1}{Introduction}
Surface transformation as a tool for studying the flow of a system of
differential equations have been introduced by Poincar\'e
(\dbiref{Poincare-1892}, Vol.~III, ch.~XXXIII) and deeply investigated by
Birkhoff~\dbiref{Birkhoff-1920}.  The corresponding method of {\corsivo
Poincar\'e sections} has become classical, and has been widely used also for
numerical explorations.  A special interesting case is that of a
periodic flow, for which the corresponding surface transformation is the
time--$T$ shift, $T$ being the period.  A natural question is whether a given
map may be represented by the time--$T$ flow of a differential system.  Positive
answers to this question exist for symplectic maps that are perturbations of
integrable ones (see,
e.g.,~\dbiref{Kuksin-1993},\dbiref{Kuksin-1994},\dbiref{Pronin-1997},\dbiref{Haro-2000}
and the references therein).

In the present paper I will reconsider the problem of giving a suitable
representation of a class of maps with a method that is somehow connected to, but
does not coincide with, the interpolation by a periodic flow.  I will rather use
the formalism of Lie transforms, which allows us to extend to maps the
techniques (e.g., calculation of normal forms) that are available for
differential equations.

The basic tools exploited here are not new: a short historical account
is given at the beginning of sect.~\secref{liemaps.2}.  However, most
practical applications of Lie series and Lie transforms methods are
related to differential equations.  Typical subjects are numerical
integration and coordinate transformations in perturbation theory.  In
the latter framework, in particular, these methods have proven to be
very effective also in investigating the convergence or the asymptotic
properties of perturbation series involving small divisors (see,
e.g.,~\dbiref{Giorgilli-1997.2}, \dbiref{Giorgilli-1999}, \dbiref{Giorgilli-2010}
and references therein) and in devising effective methods for
perturbation expansions via algebraic manipulation on computers (see,
e.g.,~\dbiref{Giorgilli-2012.1} and references therein).

It is a well known fact, however, that transporting the analytical
methods of normal form theory from differential equations to maps is
not straightforward.  A common remark is that the case of
differential equations is easier to deal with, which justifies the
attempts to interpolate a map with the Poincar\'e section of a flow.

A natural question is whether one can write a map in such a form
that transporting the normal form methods that work fine for flows
becomes straightforward.  Answering this question in general is a
major task, of course.  However, if one considers a class of maps
which are perturbation of integrable ones, then the question can be
positively answered.

One may consider as a basic example the Schr\"oder--Siegel problem of
iteration of analytic functions, extending it to many dimensions.
However, the method developed here applies also to other cases, e.g.,
perturbations of integrable symplectic maps as considered in the
papers quoted at the beginning of this section.  As interesting models
one may consider: the quadratic map investigated by
H\'enon~\dbiref{Henon-1969} and its generalization in higher
dimension; the twist map of an invariant annulus investigated by
Poincar\'e~\dbiref{Poincare-1912}, Birkhoff~\dbiref{Birkhoff-1913} and
Moser~\dbiref{Moser-1962}; the well know standard map.

In all these cases the map may be represented as a composition of two
maps: an integrable one (e.g., a linear one) and a near the identity
perturbation.  This is a trivial well known fact, of course.  However,
representing the integrable map as a Lie series and the perturbation
as a Lie transform allows one to implement the normal form theory as a
straightforward extension of the methods used for differential
equations, provided a suitable formula for the composition of Lie
transforms is available.  This is what I'm going to illustrate.

The paper is organized as follows.  In sect.~\secref{liemaps.2} a
short account of the methods based on Lie series and Lie transforms is
given, including the representation of a near the identity map and a
composition formula.  In sect.~\secref{liemaps.3} the main proposition
on the representation of a perturbation of an integrable map is
proven.  In sect.~\secref{liemaps.4} it is shown how to construct a
normal form algorithm for the map.  The actual construction is worked
out for two typical examples.  A technical appendix follows.

The theory is developed at a formal level.  Some hints on quantitative
applications to particular models are included at the end of
sect.~\secref{liemaps.4}.

\section{liemaps.2}{Basic tools}
The concept of Lie serie goes back to Sopus Lie.  The use of Lie
series in various problems has been widely investigated by Gr\"obner
in a series of papers after 1957~\dbiref{Groebner-1957}.  An accurate
exposition with particular emphasis on applications to numerical
integration can be found in Gr\"obner's book~\dbiref{Groebner-1960}.
The starting point, already found in Newton's work, is to express the
solution of an holomorphic system of differential equations as a
series expansion in time.  The basic idea is to use the representation
via power series as a one parameter near the identity map written in
an explicit and useful form.

Lie transform may be considered as a generalization of Lie series, in
a sense that will be clarified later.  Its usefulness as a tool in
perturbation theory has been emphasized by Hori~\dbiref{Hori-1966} and
Deprit~\dbiref{Deprit-1969}, who however paid attention in particular
to Hamiltonian systems. The underlying idea of Deprit's work is to
generate a one parameter family of near the identity coordinate
transformation using the flow of a {\corsivo non autonomous} system.
It is quite common in the milieu of Celestial Mechanics to call
{\corsivo Lie Transform} the algorithms proposed by Hori and Deprit,
reserving the name {\corsivo Lie series} to the algorithm based on the
flow of an {\corsivo autonomous} system; thus I follow the tradition.

Actually, several explicit algoritmhs for Lie transform have been
proposed by many authors.  A list of references may be found in
Henrard's paper~\dbiref{Henrard-1974}.  A similar algorithm for
Hamiltonian systems has been introduced on a purely algebraic basis
in~\dbiref{Giorgilli-1978}.  In this paper I will follow the latter
exposition, reformulating it for vector fields.  A rigorous treatment
may be found in~\dbiref{Giorgilli-1985}, which has been prompted by
the need of extending the contents of Gr\"obner's book to a more
general context.

I recall here the definitions and the properties of Lie series and Lie
transform working at a formal level, including what is needed
in order to develop Lie methods for maps.  I will omit most of the
proofs, that can be found elsewhere.  Furthermore, I will forget about
the origin of Lie series as solutions of a system of differential
equations, thus paying particular attention to the algebraic aspect.

\subsection{liemaps.2.1}{Lie derivatives}
Let $\Dscr\subset\complessi^n$ be an open domain endowed with
coordinates $x=(x_1,\ldots,x_n)$, and let $X(x)=(X_1,\ldots,X_n)$ be a
complex holomorphic vector field on $\Dscr$.  Let $\flusso_X^t$ denote
the time--$t$ flow generated by $X$.

The Lie derivative of a holomorphic function $f(x)$ at the point $x$
under the flow $\flusso_X^t$ is the new function
$$
\lie{X}f = \der{}{t}\bigl(\flusso_X^t f\bigr) \Big|_{t=0}
\formula{lie.1}
$$
where $\bigl(\flusso_X^t f\bigr)(x) =\bigl(f\circ\flusso_X^t\bigr)(x)$.

Similarly, the Lie derivative of a holomorphic vector field $v$ is
the new vector field
$$
\liec{X} v = \der{}{t} \bigl(\flusso_X^t v\bigr) \Big|_{t=0}
\formula{lie.2}
$$
where $\bigl(\flusso_X^t v\bigr)(x)
=\bigl(\diff \flusso_X^t\bigr)^{-1}\bigl(v\circ\flusso_X^t\bigr)(x)$.
It is well known that one has
$$
\liec{X} v = \poisson{X}{v}\ ,
\formula{lie.3}
$$
where $\poisson{X}{v}$ is the commutator between the vector fields 
$X$ and $v$.  

The Lie derivative is a linear operator mapping the space of
holomorphic functions (respectively holomorphic vector fields) into
itself.  It is also immediate to check also that the property
$$
\lie{\alpha X+\beta Y}=\alpha\lie{X}+\beta\lie{Y}
$$
holds true, where $X,\,Y$ are vector fields and $\alpha,\,\beta$ are
complex numbers.  Further useful properties are the following.  For
two functions $f,\,g$ the Leibniz rule applies, namely
$$
\lie{X}(fg)=f\lie{X}g+g\lie{X}f\ ,\quad
\lie{X}^s(fg) = \sum_{j=0}^{s} {{{s}\choose{j}}}
\bigl(\lie{X}^j f\bigr)\,\bigl(\lie{X}^{s-j} g\bigr)
\ ,\quad s\ge 1\ .
$$
For two vector fields $v,\,w$ one has
$$
\liec{X}\{v,w\} = \bigl\{\liec{X} v, w\bigr\} +
                 \bigl\{v, \liec{X} w\bigr\}\ .
$$
Finally, denoting by 
$[\lie{X},\lie{Y}] = \lie{X}\lie{Y} - \lie{Y}\lie{X}$ the commutator
between the Lie derivatives with respect to the vector fieds $X,\,Y$
one has
$$
[\lie{X},\lie{Y}] = \lie{\{X,Y\}}\ .
$$
The latter two properties are just different writings of Jacobi's
identity for the commutator between vector fields, namely
$\{X,\{v,w\}\}+\{v,\{w,X\}\}+\{w,\{X,v\}\}=0$.  In particular, if the
vector fields $X,\,Y$ do commute, namely if $\{X,Y\}=0$ then we have
$[\lie{X},\lie{Y}]=0$.

It will also be useful to write the explicit expression of the Lie
derivatives in coordinates.  For a function $f$ one has
$$
\lie{X} f=\sum_{j=1}^{n} X_j\derpar{}{x_j} f\ .
\formula{diffeq.6}
$$
For a vector field one gets the expression of the commutator, namely
$$
 \bigl(\lie{X}v\bigr)_j = \sum_{l=1}^{n} \left(X_l\derpar{v_j}{x_l}
     - v_l\derpar{X_j}{x_l}\right)\ ,
$$
where the l.h.s.~is the $j$--th component of the vector field
$\lie{X}v$.

\subsection{liemaps.2.2}{Lie series and Lie transform}
Let again $\Dscr\subset\complessi^n$ be an open domain, and let
$X=(X_1,\ldots,X_n)$ be a holomorphic vector field.  The Lie series
operator is defined as 
$$
\exp(\lie{X}) = \sum_{s\ge 0}\frac{1}{s!} \lie{X}^s
\formula{lser.1}
$$
A family of near the identity transformations depending on a parameter
$\epsilon$ may be constructed as
$$
y = \exp\bigl(\epsilon\lie{X}\bigr) x
\formula{lser.2}
$$
or, in explicit form for the coordinates,
$$
y_{j} =  \exp(\epsilon\lie{X}) x_j = x_j 
 + \epsilon X_j(x) + \frac{\epsilon^2}{2}\lie{X} X_j(x) +\ldots
\ ,\quad j=1,\ldots,n\ .
\formula{lser.3}
$$

The Lie transform is introduced as follows.  Let
$X=\{X_1,X_2,\ldots\}$ be a sequence of holomorphic vector fields,
that I will call the {\corsivo generating sequence}; here the lower
index labels the element of the sequence, not the component of the
field in coordinates, that will be denoted, e.g., by $X_{1,j}$.  The
Lie transform operator is defined as
$$
T_{X} = \sum_{s\ge 0} E^{X}_s\ ,
\formula{trslie.1}
$$
where the sequence $E^{X}_s$ of linear operators in recursively
defined as
$$
E^{X}_0 = \uno\ ,\quad 
E^{X}_s = \sum_{j=1}^{s} \frac{j}{s}\lie{X_j} E^{X}_{s-j}\ .
\formula{trslie.2}
$$
The superscript in $E^{X}$ is introduced in order to specify which
sequence of vector fields is intended.  However, I will remove it when
unnecessary.  By letting the sequence to have only one vector field
different from zero, e.g., $X=\{0,\ldots,0,X_k,0,\ldots\}$ it is
easily seen that one gets $T_X=\exp\bigl(\lie{X_k}\bigr)$.

Writing $\epsilon^s X_s$ in place of $X_s$ a one parameter family of
near the identity transformations may be defined as
$$
y = T_X x\ ,
\formula{trslie.5}
$$
i.e., in coordinates,
$$
y_j = x_j + \epsilon X_{1,j}(x) 
 + \epsilon^2 \left[\frac{1}{2} \lie{X_1} X_{1,j}(x) + X_{2,j}(x)\right] 
  +\ldots
\ ,\quad j=1,\ldots,n\ .
\formula{trslie.6}
$$
Here I used the $\epsilon$ expansion in order to make clear the
connection with classical methods based on power expansion in a small
parameter.  However, in many cases it is convenient to just consider
the vector field $X_s$ to be ``small of order $s$'' in some
appropriate sense (e.g., using a norm) the order being implicit in the
label of the field.  Adding the powers of $\epsilon$ the reader will
easily check that every term in the definition of the operator
$E^{X}_s$ carries a factor $\epsilon^s$, so that $E^{X}_s$ is of order
$s$.  An equivalent precedure is to determine the order of $E^{X}_s$
as the sum of the indices of $\lie{X_j} E^{X}_{s-j}$, which is $s$
indeed.  In the rest of the paper I will remove the parameter, unless
it has a particular meaning.  If the reader gets confused, he or she
may just rewrite a formula by adding the powers of $\epsilon$, check
that everyting is put in the correct order, and then set $\epsilon=1$.

The Lie series and Lie transform are linear operators acting on the
space of holomorphic functions and of holomorphic vector fields on the
domain $\Dscr$.  They preserve products between functions and
commutators between vector fields, i.e., if $f,\,g$ are functions and
$v,\,w$ are vector fields then one has
$$
T_X (fg)
=T_X f\cdot T_X g\ ,\quad
T_X \{v,w\} = \bigl\{T_X v, T_X w\bigr\}\ .
\formula{lser.6}
$$
Here, replacing $T_X$ with $\exp\bigl(\lie{X}\bigr)$ gives the corresponding
property for Lie series.  Moreover both operators are invertible.  The
inverse of $\exp\bigl(\lie{X}\bigr)$ is $\exp\bigl(-\lie{X}\bigr)$,
which is a natural fact if one recalls the origin of Lie series as a
solution of an autonomous system of differential equations.   The
inverse of $T_X$ takes the form
$$
\eqalign{
\bigl(T_{X}\bigr)^{-1} = 
&\sum_{s\ge 0} G^{X}_j\ ,
\cr
& G^{X}_0 = {\mat 1}\ ,\quad 
  G^{X}_s = -\sum_{j=1}^{s} \frac{j}{s} G^{X}_{s-j} \lie{X_j}\ . 
\cr
}
\formula{trslie.3}
$$

I come now to a remarkable property which justifies the usefulness of
Lie methods in perturbation theory.  I will adopt the name {\corsivo
exchange theorem} introduced by Gr\"obner.  Let $f$ be a function and
$v$ be a vector field. {\corsivo Consider the near the identity
transformation~\frmref{trslie.5} (or~\frmref{lser.2} for Lie series)
and denote by $\Jmat$ the differential of the
transformation~\frmref{trslie.5} (or~\frmref{lser.2}), namely, in
coordinates, the jacobian matrix with elements $J_{j,k}
=\derpar{y_{j}}{x_k}$.  Then one has
$$
f(y)\big|_{y=T_X x} = 
 \bigl(T_X f\bigr) (x)\ ,\quad
\Jmat^{-1} v(y)\Big|_{y= T_X x} = 
 \bigl(T_X v\bigr)(x)\ .
\formula{trslie.4}
$$
}

\noindent
This result should be interpreted as follows.  On the l.h.s.~of the
equalities there is the transformed function (resp.~vector field)
calculated via the usual method of substitution of variables.  The
r.h.s.~is the transformed function (resp.~vector field) via the Lie
transform, where the variables are renamed as $x$.  The claim is that
both operations give the same result.  The remarkable fact is that if
one uses expansions order by order, then the r.h.s.~gives it in a
straightforward way for the transformed function (resp.~vector field)
in terms of Lie derivatives, thus requiring only operations that are
easily performed, e.g., via algebraic manipulation on computers.
Obtaining the same result via substitution of variables is a
definitely longer process, unless one stops the expansion at very low
order.  This claim may appear a little obscure, but it will turn out
to be immediately evident if one writes the expansions as power series
in a parameter.

The proof of the identities~\frmref{trslie.4} may be worked out using
the algebraic properties of $T_X$.  For, a holomorphic function (or
vector field) may be expanded in power series of the variables, and by
exploiting the linearity and the preservation of product the operator
$T_X$ can be moved from the variables (the substitution) to the whole
function.  This justifies the name ``exchange theorem'', since the
symbol of the Lie transform operator is exchanged with the symbol of
the function.

\subsection{liemaps.2.3}{Representation of a near the identity transformation}
Let now a near the identity transformation be given in the form
$$
y_j = x_j + \epsilon\phi_{1,j}(x) + \epsilon^2\phi_{2,j}(x) + \ldots
\ ,\quad j=1,\ldots,n\ ,
\formula{lser.7}
$$
where $\phi_{1,j}(x),\,\phi_{2,j}(x),\ldots$ are holomorphic
functions.

\proposition{trslie.10}{Any transformation of the form~\frmref{lser.7}
may be written as a Lie transform of the coordinates $y = T_{X}x$ with
the generating sequence
$$
X_{1,j} = \phi_{1,j}\ ,\quad
X_{r,j} = \phi_{r,j} - 
  \sum_{k=1}^{r-1} \frac{k}{r} \lie{X_k} E_{r-k} x_j\ ,\quad
j=1,\ldots,n\,,\>r\gt 2\ .
\formula{trslie.7}
$$
}\endclaim

\noindent
In view of this proposition it may appear that Lie transform is more
general and more attractive than Lie series, since the latter can not
represent any near the identity transformation.  However, a similar
result may be obtained if one makes use of the {\corsivo composition
of Lie series}.

Let again $\{X_1,\,X_2,\ldots\}$ be a sequence of holomorphic vector
fields.  Consider the sequence of transformations
$\{S^{(0)}_X,\,S^{(1)}_X,\,S^{(2)}_X,\ldots\}$ recursively defined as
$$
S^{(0)}_X = 1\ ,\quad
S^{(r)}_X = \exp\bigl(\lie{X_r}\bigr) \circ S^{(r-1)}_X\ .
\formula{lser.8r}
$$
We may well consider in formal sense the limit
$$
S_X = \ldots\circ\exp\bigl(\lie{X_r}\bigr)\circ\ldots\circ
     \exp\bigl(\lie{X_2}\bigr)\circ
      \exp\bigl(\lie{X_1}\bigr)
\formula{lser.8}
$$
as an operator obtained by composition of Lie series.  A noticeable
fact, pointed out by Fass\`o~\dbiref{Fasso-1989}, is that one has
$$
S_X = \sum_{s\ge 0} :E^{X}_s:
\formula{lser.9}
$$
where $:E^{X}_s:$ denotes a reordering of $E^{X}_s$ in the following sense:
a composition of $r$ Lie derivatives
$\lie{X_{k_1}}\circ\lie{X_{k_2}}\circ\ldots\circ\lie{X_{k_r}}$ is
reordered as 
$$
:\lie{X_{k_1}}\circ\lie{X_{k_2}}\circ\ldots\circ\lie{k_r}:{\ } 
 =\> \lie{X_{\sigma_1(k)}}\circ\lie{X_{\sigma_2(k)}}\circ\ldots\circ\lie{X_{\sigma_r(k)}}
$$
where $\sigma(k)=\{\sigma_1(k),\ldots,\sigma_r(k)\}$ is any
permutation of $k=\{k_1,\ldots,k_r\}$ such that
$\sigma_1(k)\ge\sigma_2(k)\ge\ldots\ge\sigma_r(k)$.  E.g.,
$:\lie{X_1}\lie{X_2}:\ =\lie{X_2}\lie{X_1}$.

\proposition{lser.10}{Any transformation of the form~\frmref{lser.7}
may be represented via a composition of Lie series of the form $y =
S_X x$ with a generating sequence $\{X_1,\,X_2,\ldots\}$ that can be
explicitely determined with a recursive procedure.  }\endclaim

\noindent
As the reader will notice, an explicit expression for the vector
fields is missing in the statement.  Such an expression may be
produced exploiting~\frmref{lser.9}, but it turns out to be quite
useless, because it has a non recursive form.  However, a recursive
procedure for determining $X_1,\,X_2,\ldots$ may be easily constructed
by trying the first steps.  Setting $X_{1,j} =\phi_{1,j}$ one has
$$
y_j-\exp\bigl(\lie{X_{1}}\bigr) x_j 
 = \phi_{2,j} - \frac{1}{2}\lie{X_1}X_{1,j} + \ldots
$$
where the dots denote terms at least of third order.  Setting $X_{2,j}
=\phi_{2,j} -\frac{1}{2}\lie{X_1}X_{1,j}$ one gets that
$y_j-\exp\bigl(\lie{X_{2}}\bigr)\circ\exp\bigl(\lie{X_{1}}\bigr) x_j$
starts with terms at least of third order, which are used in order to
determine $X_3$, and so on.  Such a procedure is easily implemented,
e.g., via computer algebra, and it is in fact also the scheme of proof
of the proposition.  I add just a remark concerning the actual use of
propositions~\proref{trslie.10} and~\proref{lser.10}.  Apparently, it
seems unreasonable to do extra work in order to calculate the
generating sequence of the transformation.  However, if one wants to
transform either a function or a vector field, pushing the calculation
at high orders, then using the exchange theorem turns out to be
definitely more effective than performing a substitution.

\subsection{liemaps.2.4}{Composition formul{\ae}}
I begin with a formula for the commutation between two Lie
transforms.  {\corsivo Let $X,\,Y$ be generating sequences.  Then one has
$$
T_X \circ T_Y = T_{W} \circ T_X\ , 
\formula{trslie.8}
$$
with the generating sequence $W=\{T_X Y_1,\,T_X Y_2,\ldots\}$.}  Here,
the Lie trasform $T_X$ may be replaced by $\exp\bigl(\lie{X}\bigr)$ in
case $X$ is a vector field, and similarly for $Y$.

Generally speaking the formula above does not seem very interesting
because the vector field $W$ turns out to be itself a series.
However, it is useful in some case, typically when one is able to put
$W$ in a manageable form by explicitly calculating the sum of the
series.  It will be used later, in sect.~\sbsref{liemaps.3.2}.

It is instead more interesting to to observe that since $y = T_X\circ
T_Y x$ is a near the identity transformation then in view of
proposition~\proref{trslie.10} there exists a generating sequence $Z$
such that $y=T_Z x$.  The following proposition gives an explicit
expression for $Z$.

\proposition{trslie.11}{Let $X,\,Y$ be generating sequences.  Then one
has $T_X\circ T_Y = T_Z$ where $Z$ is the generating sequence
recursively defined as
$$
Z_1 = X_1 + Y_1\ ,\quad
Z_s = X_s + Y_s + \sum_{j=1}^{s-1} \frac{j}{s} E^{X}_{s-j} Y_j\ .
\formula{trslie.12}
$$
}\endclaim

\noindent
The proof requires a long sequence of tedious calculations that can
hardly be found in previous papers. Thus I include it in
appendix~\appref{liemaps.a}.  Composition formul{\ae} for any
combination of Lie series and Lie transforms are easily obtained by
suitably elaborating formula~\frmref{trslie.12}.  E.g., if $X,\,Y$ are
vector fields then one has
$\exp\bigl(\lie{X}\bigr)\circ\exp\bigl(\lie{Y}\bigr)= T_{W}\,$, where
$W=\{W_1,\,W_2,\ldots\}$ is the generating sequence
$$
W_1 = X + Y\ ,\quad
W_s = \frac{1}{s!} \lie{X}^{s-1} Y\ .
\formula{trslie.13a}
$$
The latter formula reminds the well known Baker--Campbell--Hausdorff
composition of exponentials.  The difference is that the result is
expressed as a Lie transform instead of an exponential.

\section{liemaps.3}{Representation of a map by a composition of Lie Transforms}
I shall consider two cases.  The first one is a map in a neighbourhood
of an equilibrium, which may be expanded in Taylor series as
$$
z' = \Lambdamat z + v_1(z) + v_2(z) + \ldots\ ,\quad z\in\complessi^n
\formula{map.1}
$$
where $\Lambdamat$ is a $n\times n$ complex matrix and $v_s(z)$ is a
homogeneous polynomial af degree $s+1$.  The second example is
a real analytic map
$$
\eqalign{
\phi' &= \phi + \omega(I) + \epsilon f_1(\phi,I) 
                           + \epsilon^2 f_2(\phi,I) +\ldots 
\cr
I' &= I + \epsilon g_1(\phi,I) + \epsilon^2 g_2(\phi,I) + \ldots
\cr
}
\formula{map.2}
$$ 
where $(\phi,I)\in\toro^n\times\Gscr$,  with $\Gscr\subset\reali^m$
open, and $\epsilon$ is a small perturbation parameter.

\subsection{liemaps.3.1}{The unperturbed map}
By unperturbed map I mean here (as usual) either the linear part
of~\frmref{map.1} or the map~\frmref{map.2} with $\epsilon=0$, which is
a Kronecker map on a family of invariant tori parameterized by the
actions $I$.

Let me start with the linear part $z' = \Lambdamat z$ of the
map~\frmref{map.1}.  Let $\Lambdamat=e^{\Amat}$ with a complex
$n\times n$ matrix $\Amat$.  For the vector field $\Amat z$ one has
$\lie{\Amat z} z = \Amat z \,, \ldots,\> \lie{\Amat z}^s z = \Amat^s
z$, and so also $\exp\bigl(\lie{\Amat z}\bigr) z =\sum_{s\ge
0}\frac{1}{s!}\Amat^s z = e^{\Amat} z$.  Using the
exchange theorem we may transform a function $f(z')$ and a vector
field $v(z')$ as
$$
\vcenter{\openup1\jot\halign{
\hfil$\displaystyle{#}$
&$\displaystyle{#}$\hfil
&$\displaystyle{#}$\hfil
\cr
\exp\bigl(\lie{\Amat z}\bigr) f(z)
&= f(z')\Big|_{z'=e^{\Amat} z}
&= f(\Lambdamat z)\ ,
\cr
\exp\bigl(\lie{\Amat z}\bigr) v(z)
&= e^{-\Amat} v(z')\Big|_{z'=\exp(\lie{\Amat z}) z}
&= \Lambdamat^{-1} v(\Lambdamat z)\ .
\cr
}}
\formula{map.1a}
$$
The simplest case occurs when $\Lambdamat
=\diag(\lambda_1,\ldots,\lambda_n)$ is a diagonal matrix.

I point out that here the exchange theorem is used only in order to
represent a transformation of a function and of a vector field as the
action of a Lie series operator.  Of course, in this case everybody
would perform the transformations via a direct substitution as in the
last member of the formul{\ae} above, since this is actually the sum
of the Lie series in closed form and requires no further expansion.

I come now to considering the map~\frmref{map.2}, that for
$\epsilon=0$ writes
$$
\phi' = \phi + \omega(I)\ ,\quad I'=I.
\formula{map.3}
$$
Introducing the $(n+m)$--dimensional vector field
$\Omega=\bigl(\omega(I), 0\bigr)$ one immediately gets
$$
\lie{\Omega} (\phi, I)  =  \bigl(\omega(I), 0\bigr)\ ,
\quad \lie{\Omega}^s \bigl(\omega(I), 0\bigr) = (0,0)\quad 
{\rm for}\ s\gt 1\>.
$$
Thus, the map may be rewritten as $(\phi',I')
=\exp\bigl(\lie{\Omega}\bigr)(\phi,I)$.  Using again the exchange
theorem as above we may transform a function $f(\phi,I)$ or a
$(n+m)$--dimensional vector field $v(\phi,I)$ as
$$
\vcenter{\openup1\jot\halign{
\hfil$\displaystyle{#}$
&$\displaystyle{#}$\hfil
&$\displaystyle{#}$\hfil
\cr
\exp\bigl(\lie{\Omega}\bigr) f(\phi,I) 
&= f(\phi',I')\Big|_{(\phi',I')=\exp(\lie{\Omega})(\phi,I)}
&= f\bigl(\phi+\omega(I),I\bigr)\ ,
\cr
\exp\bigl(\lie{\Omega}\bigr) v(\phi,I)
&= \Jmat^{-1} v(\phi',I')\Big|_{(\phi',I')=\exp(\lie{\Omega})(\phi,I)}
&= \Jmat^{-1} v\bigl(\phi+\omega(I),I\bigr)\ ,
\cr
}}
\formula{map.2a}
$$
where $\Jmat$ and its inverse $\Jmat^{-1}$ are the jacobian
block matrices
$$
\Jmat     = \pmatrix{\Imat_n & \Bmat\phantom{_{m}} \cr
                      0      & \Imat_m}\ ,\quad
\Jmat^{-1} = \pmatrix{\Imat_n & -\Bmat\phantom{_{m}} \cr
                      0      &  \phantom{-}\Imat_m}\ ,\quad
\Bmat = \biggl\{\derpar{\omega_j}{I_l}\biggr\}_{1\le j\le n,\,1\le l\le m}\ ,
$$
$\Imat_s$ denoting the $s\times s$ identity matrix.

\subsection{liemaps.3.2}{The full map}
The maps~\frmref{map.1} or~\frmref{map.2} considered in the last
section are just interesting examples.  In more general terms one can
consider an unperturbed map $x'=f_0(x)$ which has some nice properties
and can be represented as a Lie transform $x'
=\exp\bigl(\lie{X}\bigr)$ with an appropriately defined vector field
$X$.  In view of the particular role of the latter map I will denote
$\Rmat =\exp\bigl(\lie{X}\bigr)$.

Then one considers a perturbed map
$$
x' = \Rmat x + f_1(x) + f_2(x) + \ldots
\formula{map.4}
$$
where $f_s(x)$ is of order $s$ in some reasonable sense.  E.g.,
$f_s(x)$ should be replaced by $v_s(z)$ in~\frmref{map.1}, and by
$\epsilon f_s(\phi,I)$ in~\frmref{map.2}, with the obvious change of
the symbols for the coordinates.  

\proposition{map.13}{Consider the map~\frmref{map.4} where $\Rmat$
is a Lie series operator.  Then there exist generating sequences of
vector fields $V(x)= \bigl\{V_1(x),V_2(x),\ldots\bigr\}$ and $W(x)
=\bigl\{W_1(x),W_2(x),\ldots\bigr\}$ with $W_s =\Rmat V_s$ such that
one has both
$$
x' = \Rmat\circ T_{V} x\quad {\rm and}\quad
x' = T_{W}\circ\Rmat\, x
\formula{map.5}
$$
}\endclaim

\proof
Using the linearity of the Lie series operator $\Rmat$ rewrite the
map~\frmref{map.5} as
$$
x' = \Rmat \bigl(x + \Rmatinv \,(f_1+f_2+\ldots)\bigr)\ .
$$
This is a representation of the map as the composition of two
operations, namely
$$
\tilde\phi(x) = x + \Rmatinv \,(f_1+f_2+\ldots) \ ,\quad
x' =\Rmat \tilde\phi(x)\ .
$$
In view of proposition~\proref{trslie.10} one may determine the
generating sequence $V(x) =\bigl\{V_1(x),V_2(x),\ldots\bigr\}$ such
that $\tilde\phi(x)=T_{V} x$.  Applying $\Rmat$ to both members and
using linearity one readily gets $x' = \Rmat \circ T_{V} x\,$, namely
the first of~\frmref{map.5}.  Using the identity~\frmref{trslie.8} of
proposition~\proref{trslie.11}, which clearly applies also to the Lie
series operator $\Rmat$, one gets $\Rmat \circ T_{V} =
T_{W} \circ \Rmat$ with $W$ as in the statement, which gives the second
of~\frmref{map.5}.\endproof

\subsection{liemaps.3.3}{Conjugating maps}
I come now to the following question. {\corsivo Let two maps
$$
x'=T_{W}\circ\Rmat\, x\ ,\quad y'=T_{Z}\circ \Rmat\, y
\formula{map.6}
$$ 
be given, where $\Rmat$ is an invertible Lie series operator and
$W=\{W_1,W_2,\ldots\}$, $Z=\{Z_1,Z_2,\ldots\}$  are generating
sequences.  To find whether the maps are conjugated by  a holomorphic
near the identity transformation}
$$
y = x + \phi_1(x) + \phi_2(x) +\ldots\ .
\formula{map.7}
$$
Using the same operator $\Rmat$ in both maps means only that the
unperturbed maps are trivially conjugated.  In view of
propositions~\proref{trslie.10} and~\proref{lser.10} it is natural to
consider the transformation~\frmref{map.7} as generated by either a
Lie transform or a Lie series.

\proposition{map.10}{Let  $X=\{X_1,X_2,\ldots\}$ be a generating
sequence of the near the identity transformation $y=T_X x$.  Then the
maps~\frmref{map.6} are conjugated if 
$$
T_{W}\circ T_{\Rmat X} = T_{X} \circ T_{Z}\ .
\formula{map.9}
$$
More explicitly, the following relations must be satisfied:
\formdef{map.20c}
\formdef{map.20b}
$$
\displaylines{
\frmref{map.20c}\hfill
\Dmat X_1 = Z_1 - W_1\ ,\quad \Dmat = \Rmat-\uno
\hfill\cr
\frmref{map.20b}\hfill
\Dmat X_s =
 Z_s - W_s + \sum_{j=1}^{s-1}\frac{j}{s} 
  \bigl(E^{X}_{s-j} Z_j - E^{W}_{s-j} \Rmat X_j\bigr)  
\ ,\quad s\gt 1\ .
\hfill\cr
}
$$}\endclaim

\proof
In the l.h.s.~of the map $y'=T_{Z}\circ\Rmat\, y$ set $y'=T_{X} x'$
and then substitute $x' =T_{W}\circ\Rmat\, x$.  This must be the same
as substituting $y=T_{X} x$ in the r.h.s.~of the map.  Thus the
identity
$$
T_{X} x'\Big|_{x' = T_{W}\circ\Rmat\, x} 
 = T_{Z}\circ\Rmat\, y\Big|_{y = T_{X} x} 
$$
must be true, and by the exchange theorem this gives
$T_{W}\circ\Rmat\circ T_{X} = T_{X}\circ T_{Z}\circ\Rmat$.  In view
of~\frmref{trslie.8} one has $\Rmat\circ T_{X} = T_{\Rmat
X}\circ\Rmat$, so that~\frmref{map.9} is readily found in view of the
invertibility of $\Rmat$.  Apply now proposition~\proref{trslie.11} to
both members of the latter equality.  By the first
of~\frmref{trslie.12} we get $W_1 + \Rmat X_1 = X_1 + Z_1$,
namely~\frmref{map.20c}.  For $s\gt 1$ we get
$$
W_s +\Rmat X_s + 
 \sum_{j=1}^{s-1}\frac{j}{s} E^{W}_{s-j} \Rmat X_j =
X_s + Z_s +
 \sum_{j=1}^{s-1}\frac{j}{s} E^{X}_{s-j} Z_j\ ,
$$
from which~\frmref{map.20b} readily follows.\endproof

A similar result for Lie series holds true if in~\frmref{map.9} one
replaces the Lie transform $T_{X}$ with the infinite composition of
Lie series $S_{X}$ as given by~\frmref{lser.8}.
However,~\frmref{map.20b} must be restated in a more elaborated
manner, proceeding step by step.  Let me say that the
maps~\frmref{map.9} are conjugated up to order $r$ in case there
exists a finite generating sequence $X=\{X_1,\ldots,X_r\}$ such that
the transformation $y=S^{(r)}x$ makes the difference between the maps
to be of order higher than $r$, i.e.,
$$
S^{(r)}_{X} x'\Big|_{x' = T_{W}\circ\Rmat\, x} 
 - T_{Z}\circ\Rmat\, y\Big|_{y = S^{(r)}_{X} x} = \Oscr(r+1)\ ,
$$
where $S^{(r)}$ is defined by~\frmref{lser.8r}.  The maps are
trivially conjugated up to order $r$ if the generating sequences
$W,\,Z$ coincide up to order $r$, i.e., $W_1=Z_1,\ldots,W_r=Z_r$.  For
this implies $T_{W} x - T_{Z} x = O(r+1)$.

\proposition{map.20}{Let the generating sequences of the
maps~\frmref{map.6} coincide up to order $r-1$ and let $X_r$ be a
vector field of order $r$ generating the near the identity
transformation $y=\exp\bigl(\lie{X_r} x\bigr)$.  Then the maps are
conjugated up to order $r$ if
$$
T_{W}\circ \exp\bigl(\lie{\Rmat X_r}\bigr) = 
   \exp\bigl(\lie{X_r}\bigr) \circ T_{Z}\ .
\formula{map.21a}
$$
More explicitly, the following relations must be satisfied:
\formdef{map.21b}
\formdef{map.21c}
$$
\vbox{
\displaylines{
\frmref{map.21b}\hfill
\Dmat X_r  = Z_r - W_r\ ,\quad \Dmat = \Rmat-\uno\ ;
\hfill\cr
\frmref{map.21c}\quad
 Z_{s} - W_{s} =
  \frac{r}{s} E^{W}_{s-r} \Rmat X_r 
   -\sum_{j=1}^{\lfloor (s-1)/r\rfloor} \frac{s-jr}{s\cdot j!} \lie{X_r}^{j} Z_{s-jr}
\qquad
{{\rm for}\ s\gt r\>.}
\hfill\cr
}}
$$}\endclaim

\noindent
The proof of~\frmref{map.21a} is a straightforward adaptation of that
proposition~\proref{map.10} using a generating sequence $X$ with all
elements zero except $X_r\,$.  For, in this case one has
$T_{X}=\exp\bigl(\lie{X_r}\bigr)$, so that $E^{X}_s=0$ if $s$ is not a
multiple of $r$.  Removing from~\frmref{map.20b} all vanishing term,
with some patience one obtains~\frmref{map.21c}.

\section{liemaps.4}{Normal form algorithm} 
A standard and useful tool in perturbation theory is the construction
of a normal form for either a map or a system of differential
equations.  This is a classical and widely investigated problem, so I
will limit the discussion to indicating how the known methods may be
revisited within the scheme of representation of maps presented in
this paper.  Constructing a normalization algorithm for maps is indeed
an easy matter in view of the results of sect.~\sbsref{liemaps.3.3}.
Actually, two different algorithms may be devised, the first one based
on Lie transform, the second one based on composition of Lie series.

\subsection{liemaps.4.1}{A general formulation}
Let me start with the Lie transform.  The key point is that
formula~\frmref{map.9} must be considered as an equation for the
generating sequence $X$, which must be so determined that the
transformed sequence $Z$ has some nice property that characterizes it
as being in normal form.  Look now at equations~\frmref{map.20b},
which are just a rewriting of~\frmref{map.9} order by order.  The
reader will immediately see that for every $s\gt 1$ one has to solve
recursively an equation of the form
$$
\Dmat X_s + Z_s = \Psi_s
\formula{nrm.1}
$$
where $\Psi_s$ is known, since it is determined by
$X_1,\ldots,X_{s-1}$ and $W_1,\ldots,W_{s-1}$, which are known.  Thus,
the problem is only that the prescription that $Z$ has a normal form
should be imposed so that eq.~\frmref{nrm.1} may be solved for $X_s$
and $Z_s$.  This is the standard problem in normal form theory for
both maps and differential equations.  If the process can be worked
out, at least formally, then the generating sequence $X$ produces a
coordinate transformation $y = T_{X} x$, with inverse $x=T_{X}^{-1}
y$, such that the map in the new coordinates writes $y' =
T_{Z}\circ\Rmat\, y$.

The algorithm based on composition of Lie series appears to be more
elaborated, since it requires using proposition~\proref{map.20} as an
iteration step.  Precisely, one constructs an infinite sequence
$\{W^{(r)}\}_{r\ge 0}$ of generating sequences, with $W^{(0)}=W$, and
a generating sequence $X=\{X_1,X_2,\ldots\}$ such that for every $r\gt
0$ the generating sequences $W^{(r-1)}$ and $W^{(r)}$ coincide up to
order $r-1$, and $Z_1=W^{(r)}_1,\ldots,Z_r=W^{(r)}_r$ are in normal
form and do not change with the next iteration.
Proposition~\proref{map.20} shows how to determine $X_r$ and
$W^{(r)}_{r}=Z_r$ by solving eq.~\frmref{map.21b} with $W^{(r-1)}_r$
in place of $W_r$.  Then the whole generating sequence $W^{(r)}$ is
constructed as given by~\frmref{map.21c}, putting $W^{(r-1)}$ in place
of $W$ and $W^{(r)}$ in place of $Z$.  The normal form is thus
determined step by step.

At first sight, the reader may think that this is a too complicated
process.  However, by implementing the algorithm using algebraic
manipulation he or she will realize that there is no substantial
increase of complexity with respect to the algorithm based on a single
Lie transform, and that in some cases the composition of Lie series
may even be more effective.

\subsection{liemaps.4.2}{Back to examples}
Let me illustrate how a normal form may be constructed for the
examples~\frmref{map.1} and~\frmref{map.2}.  Actually, this means that
I should explain how to characterize the normal form and how to solve
eq.~\frmref{nrm.1}.  

Assume that the matrix $\Lambdamat$ in~\frmref{map.1} has a diagonal
form, namely $\Lambdamat=\diag(\lambda_1,\ldots,\lambda_n)$, and
define $\lambda_j =e^{\mu_j +i\omega_j}$.  This is a $n$--dimensional
version of the problem of iteration of analytic maps investigated by
Schr\"oder~\dbiref{Schroeder-1871}, who gave the formal solution for
the case $n=1$.  Clearly one has
$$
\Rmat = \exp\bigl(\lie{\Omegamat x}\bigr)\ ,\quad
 \Omegamat = \diag(\mu_1 +i\omega_1,\ldots,\mu_n +i\omega_n)\ ,
$$
and applying the method of proposition~\proref{map.13} the generating
sequence $W$ may be determined so that $W_s$ is a homogeneous
polynomial of degree $s+1$.  Following Schr\"oder one tries to
conjugate the map to its linear part.  This means that the normal form
of the map should be $\zeta'=\Lambdamat\zeta$, which in the Lie
transform representation means that one wants $Z=\{0,0,\ldots\}$, the
null sequence.  Thus, according to~\frmref{nrm.1}, the generating
sequence is determined by solving for $X_s$ the equation
$$
\Dmat X_s = \Psi_s\ ,\quad
 \Dmat = \exp\bigl(\lie{\Omegamat x}\bigr) - \uno\ , 
\formula{nrm.2}
$$
where $\Psi_s(x_1,\ldots,x_n)$ is a homogeneous polynomial of degree
$s+1$.  The relevant property is that $\Dmat$ is diagonal on the basis
of monomials $x^k\evet_j = x_1^{k_1}\cdot\ldots\cdot
x_n^{k_n}\evet_j$, where $(\evet_1,\ldots,\evet_n)$ is the canonical
basis of $\complessi^n$.  For, in view of the second of~\frmref{map.1a}
one has
$$
\Dmat\, x^k\evet_j 
 = \bigl(e^{\langle k,\mu+i\omega\rangle - \mu_j-i\omega_j}-1\bigr)\,
    x^k\evet_j \ .
$$
Thus determining the vector field $X_s$ is an easy matter if none of
the eigenvalues of $\Dmat$ is zero.  For, writing the homogeneous
polynomial vector field as
$$
\Psi_{s} = 
 \sum_{j=1}^{n}\evet_j\sum_{|k|=s+1}\psi_{j,k} x^k
$$
the solution of~\frmref{nrm.2} is readily found to be
$$
X_{s} = \sum_{j=1}^{n} \evet_j \sum_{k} 
 \frac{\psi_{j,k}}{e^{\langle k,\mu+i\omega\rangle-\mu_j-i\omega_j}-1} 
   x^k\ .
$$
Thus the map may be formally linearized if the {\corsivo nonresonance
condition}
$$
e^{\langle k,\mu+i\omega\rangle - \mu_j - i\omega_j} \ne 1
\ {\rm for}\ k\in\interi_{+}^n\>, |k|\gt 1\ {\rm and}\ j=1,\ldots,n
\formula{nrm.3}
$$
is fulfilled.  In the case $n=1$ considered by Schr\"oder a resonance
may occur only if $\lambda$ is a root of the unity.

If the resonance condition is violated then a different definition of
normal form must be introduced.  Precisely, the space $\Pscr^s$ of
homogeneous polynomials of any degree $s$ splits into two
complementary subspaces
$$
\Nscr^s = \Dmat^{-1}\{0\}\ ,\quad \Rscr^s=\Dmat(\Pscr^s)\ ,
$$
namely the kernel of $\Dmat$ and the image of $\Pscr^s$ through
$\Dmat$.  For, $\Dmat$ maps $\Pscr^s$ onto itself, so that both
$\Nscr^s$ and $\Rscr^s$ are subspaces of $\Pscr^s$, and $\Dmat$ is
diagonal.  Then the operator $\Dmat$ may be uniquely inverted on
$\Rscr^s$.  Thus eq.~\frmref{nrm.1} may be solved by splitting
$\Psi_s=\Psi^{\Nscr}_s+\Psi^{\Rscr}_s$, with obvious meaning of the
superscripts, and setting
$$
Z_s = \Psi^{\Nscr}_s\ ,\qquad 
 X_s = \Dmat^{-1}\Psi^{\Rscr}_s\ ,\quad X_s\in\Rscr^{s}\ .
$$
The latter condition makes the solution unique.  The procedure thus
described is a standard one in normal form theory.  Different
solutions may be considered, of course, depending on what one is
looking for.

The convergence of the tranformation to normal form in the non
resonant case has been proved by Siegel~\dbiref{Siegel-1942} for the
case $n=1$ under the additional hypothesis that $\lambda$ satisfies a
diophantine condition.  The work of Siegel represents a milestone for
the problem of convergencence of perturbation series with small
divisors.

Let me now come to the model~\frmref{map.2}.  In view of the
particular form of the map it is convenient to represent the vector
fields separating, so to say, the $\phi$ component from the $I$
component by writing $\pmatrix {X\cr Y\cr}$ in place of $X$,
where $X(\phi,I)$ and $Y(\phi,I)$ are a $n$--dimensional and a
$m$--dimensional vector function, respectively.  
The explicit form of the commutator is written as
$$
\left\{\pmatrix {X \cr Y\cr}, \pmatrix {V \cr W\cr}\right\} = 
\pmatrix{
\sum_{l=1}^{n} 
 \left(X_l\derpar{V}{\phi_l} - V_l\derpar{X}{\phi_l}\right) +
  \sum_{l=1}^{m} 
   \left(Y_l\derpar{V}{I_l} - W_l\derpar{X}{I_l}\right)
\cr
\sum_{l=1}^{n} 
 \left(X_l\derpar{W}{\phi_l} - V_l\derpar{Y}{\phi_l}\right) +
  \sum_{l=1}^{m} 
   \left(Y_l\derpar{W}{I_l} - W_l\derpar{Y}{I_l}\right)
\cr
}
\formula{map.22}
$$
Recall also that in this case whith the notation above one has
$$
\Dmat =  \exp\bigl(\lie{\Omega}\bigr) - \uno\ ,\quad
 \Omega = \pmatrix{\omega(I) \cr 0\cr}\ .
$$
In the r.h.s.~of~\frmref{nrm.1} we may expand
$\Psi_s(\phi,I)$ in Fourier series
as
$$
\Psi_s = \pmatrix{
 \sum_{k\in\interi^m} \alpha_{k}(I) e^{i\langle k,\phi\rangle} \cr
  \sum_{k\in\interi^m} \beta_{k}(I) e^{i\langle k,\phi\rangle} \cr
}
$$
with known coefficients $\alpha_{k}(I)$ and $\beta_{k}(I)$.  Using a
similar expansion for
$$
X_s = \pmatrix{
 \sum_{k\in\interi^m} c_{k}(I) e^{i\langle k,\phi\rangle} \cr
  \sum_{k\in\interi^m} d_{k}(I) e^{i\langle k,\phi\rangle} \cr
}\ ,
$$
in view of~\frmref{map.2a} the action of the operator $\Dmat$ is given
by
$$
\Dmat X_s = \pmatrix{
 \sum_{k\in\interi^n} 
  \bigl(e^{i\langle k,\omega(I)\rangle}-1\bigr)
   c_{k}(I) e^{i\langle k,\phi\rangle}
    - \sum_{k\in\interi^n} e^{i\langle k,\omega(I)\rangle}
     \Bmat d_{k}(I) e^{i\langle k,\phi\rangle}
\cr
\sum_{k\in\interi^n} \bigl(e^{i\langle k,\omega(I)\rangle}-1\bigr)
 d_{k}(I) e^{i\langle k,\phi\rangle}
\cr
}\ .
$$
Thus, one would be tempted to solve eq.~\frmref{nrm.1} by setting
$$
Z_s = \pmatrix{\alpha_{0}(I) \cr \beta_{0}(I)}\ ,\quad
 X_s = \pmatrix{
\sum_{k\in\interi^n\setminus\{0\}}
 \left[\frac{\alpha_{k}(I)}{e^{i\langle k,\omega(I)\rangle}-1}
  + \frac{e^{i\langle k,\omega(I)\rangle}\Bmat\beta_{k}(I)}{(e^{i\langle k,\omega(I)\rangle}-1)^2}
   \right] e^{i\langle k,\phi\rangle}
\cr
\sum_{k\in\interi^n\setminus\{0\}}
 \frac{\beta_{k}(I)}{e^{i\langle k,\omega(I)\rangle}-1}
  e^{i\langle k,\phi\rangle}
\cr 
}\ .
\formula{map.23}
$$
However, some denominator could vanish at some point $I\in\Gscr$, the
action's domain, or at least become very small.  This is indeed the
classical problem of small divisors in Celestial Mechanics, which was
well known to, e.g., Lagrange and Laplace.

Assume for a moment that no divisor actually vanishes.  This is true,
e.g., if $\omega\in\reali^n$ is a constant vector and the non
resonance condition $e^{i\langle k,\omega\rangle}\ne 1$ for
$k\in\interi^m\setminus \{0\}$ is satisfied.  Then the generating
sequence $Z(I)$ is independent of the angle variables $\phi$.
Thus one has
$$
\lie{Z_s(I)}\phi = Z_s(I)\ ,\quad
 \lie{Z_s}^r \phi = 0\ {\rm for}\ r\gt 1\ ,
$$
and the map in normal form is written as
$$
\phi' = \phi + \omega'(I)\ ,\quad
 I' = I + g'(I)
$$
where $\omega'(I)$ and $g'(I)$ are determined via the normalization
process.  Thus the dynamics of the actions $I$ is separated from that
of the angles $\phi$.

A simpler form of the normalized map is found in case the map
possesses some interesting symmetries.  Let me give an example.  Say
that a vector field ${\scriptscriptstyle\pmatrix{X\cr Y\cr}}$ is of
type $(+,-)$ or of type $(-,+)$, respectively, if it satisfies
$$
\pmatrix{X(-\phi,I)\cr Y(-\phi,I)} =
 \pmatrix{X(\phi,I)\cr -Y(\phi,I)}
\quad{\rm or}\quad
\pmatrix{X(-\phi,I)\cr Y(-\phi,I)} =
 \pmatrix{-X(\phi,I)\cr Y(\phi,I)}\ .
$$ 
With a little patience, looking at the explicit
expression~\frmref{map.22} of the commutator, one checks that the
commutators obey the rules symbolically expressed by the table
$$
\vcenter{\tabskip=0pt
\def\tablerule{\noalign{\hrule}}
\halign{
 \hbox to 4 em{\hfil$\displaystyle{#}\hfil$}\vrule\vrule
&\hbox to 4.5 em{\hfil$\displaystyle{#}\hfil$}\bigg\vert
&\hbox to 4.5 em{\hfil$\displaystyle{#}\hfil$}\vrule\vrule
\cr
\{{\cdot},{\cdot}\}
 & (+,-)
 & (-,+)
\cr
\tablerule
\tablerule
(+,-)
 & (-,+)
 & (+,-)
\cr
\tablerule
(-,+)
 & (+,-)
 & (-,+)
\cr
\tablerule
\tablerule
}}\ .
$$
Assume now that the map~\frmref{map.2} satisfies the symmetry
$$
f_s(-\phi,I) = f_s(\phi,I)\ ,\quad
 g_s(-\phi,I) = -g_s(\phi,I)\ .
$$
This means that the system is reversible.  It is not difficult to
check that in this case the generating sequence $W(\phi,I)$ is of type
$(+,-)$.  Then, with a little more patience and using induction, one
checks also that at every step of the normalization procedure one has
that $\Psi_s$ is of type $(+,-)$, and so by solving eq.~\frmref{nrm.1}
one gets $Z_s$ of type $(+,-)$ and $X_s$ of type $(-,+)$.  This
implies that in~\frmref{map.23} one has
$Z_s={\scriptscriptstyle{\pmatrix{\alpha_{0}(I)\cr 0\cr}}}$.  In turn
this implies that the normalized map is written as
$$
\phi' = \phi + \omega'(I)\ ,\quad
 I' = I\ ,
$$
representing a Kronecker map on a family of invariant tori with angles
$\omega'(I)$.

All this is formal, of course.  Making a rigorous statement, as is
known, is a definitely more complicated matter.  If the angles
$\omega(I)$ of the unperturbed map do depend on the actions $I$ then
the formal construction above is expected to fail due to the presence
of zero divisors, and an analog of Poincar\'e's theorem on non
existence of holomorphic first integrals for Hamiltonian systems of
differential equations applies.  However one can
prove, possibly with some extra condition on $\omega(I)$, that KAM
theory applies, thus showing the existence of a big set on invariant
tori carrying a Kronecker map with strongly non resonant angles.  On
the other hand, it should also be possible to prove a theorem of
Nekhoroshev's type on exponential stability.  If the angles $\omega$
are constant, then some normal form may be constructed, possibly
taking into account the resonances, but the series so constructed are
expected to be divergent, unless one looks for a Kolmogorov's normal
form on an invariant torus.  However, all this matter goes behind the
limits of the present note, which deals only with formal aspects.

\appendix{liemaps.a}{Proof of proposition~\proref{trslie.11}}
Let me state a preliminary identity.  If $X$ is a generating sequence
and $V$ a vector field then one has
$$
E^{X}_s \lie{V} 
 = \sum_{j=0}^{s} \lie{E^{X}_j V} E^{X}_{s-j}\ .
\formula{trslie.5b}
$$
The proof is worked out by induction, since the equality is trivial
for $s=0$ and moreover for $s=1$ it is just Jacobi's identity for
commutators.  Here is the complete calculation for $s\gt 1$.  The
notation is made simpler by writing $E_{j}$ in place of $E^{X}_j$,
since there is no confusion.
$$
\leqalignno{
E_s\lie{V}
&= \sum_{m=1}^{s}\frac{m}{s}\lie{X_{m}} E_{s-m}
    \lie{V}
\cr
&= \sum_{m=1}^{s}\frac{m}{s} \lie{X_{m}} 
    \sum_{j=0}^{s-m} \lie{E_j V} E_{s-m-j} 
\cr
&= \sum_{m=1}^{s} \sum_{j=0}^{s-m} \frac{m}{s}
    \Bigl(\lie{\lie{X_{m}}E_{j}V}+
     \lie{E_j V}{\lie{X_{m}} \Bigr) E_{s-m-j}}
\cr
&= \sum_{m=1}^{s} \sum_{j=0}^{s-m} \frac{m}{s}
    \Bigl(\lie{\lie{X_{m}}E_{s-j-m}V} E_j +
     \lie{E_j V}{\lie{X_{m}}E_{s-j-m}} \Bigr) 
\cr
&= \sum_{j=0}^{s-1} \frac{s-j}{s}
    \sum_{m=1}^{s-j} \frac{m}{s-j}
     \Bigl(\lie{\lie{X_{m}}E_{s-j-m}V} E_j +
      \lie{E_j V}{\lie{X_{m}}E_{s-j-m}} \Bigr) 
\cr
&=\sum_{j=0}^{s-1} \frac{s-j}{s}
   \Bigl(\lie{E_{s-j}V} E_j +
    \lie{E_j V} E_{s-j} \bigr) 
\cr
&= \sum_{j=1}^{s} \frac{j}{s} \lie{E_j V} E_{s-j} +
    \sum_{j=0}^{s-1} \frac{s-j}{s}\lie{E_j V} E_{s-j} 
= \sum_{j=0}^{s} \lie{E_j V} E_{s-j}\ .
\cr
}
$$
Jacobi's identity for commutators written as
$\lie{X}\lie{w}-\lie{w}\lie{X} = \lie{\lie{X}w}$ is used in order to
obtain the third equality.

Coming to the proof of~\frmref{trslie.12} of
proposition~\proref{trslie.11}, by definition of Lie transform one has
$$
T_{X}\circ T_{Y} = 
 \biggl(\sum_{l\ge 0} E^{(X)}_l\biggr)
  \biggl(\sum_{k\ge 0} E^{(Y)}_k\biggr)
= \sum_{s\ge 0} \sum_{m=0}^{s} E^{(X)}_{m\phantom{-}} E^{(Y)}_{s-m}\ .
\formula{trslie.13}
$$
On the other hand for the generating sequence $Z$ defined as
in~\frmref{trslie.12} one has
$$
E^{Z}_s=\sum_{l=1}^{s}{l\over s} L_{X_{l}+Y_{l}}E^{Z}_{s-l}+
 \sum_{l=2}^{s} \fraz{l}{s}
  \sum_{m=1}^{l-1} {m\over l} L_{E^{X}_{l-m}Y_{m}} E^{Z}_{s-l}\ .
\formula{trslie.8a}
$$
Thus it is enough to check that 
$$
E^{Z}_s= \sum_{m=0}^{s} E^{X}_{m\phantom{-}} E^{Y}_{s-m}\ ,\quad s\ge 0\>.
\formula{trslie.8b}
$$
I proceed by induction. For $s=0,\,1$ the equality is true.  For $s\gt
1$ calculate
$$
\incolonna{
\sum_{m=0}^{s} E^{X}_{m\phantom{-}} E^{Y}_{s-m} 
&= \sum_{m=0}^{s-1} {s-m\over s} 
    \left(E^{X}_{s-m} E^{Y}_{m\phantom{-}} 
     + E^{X}_{m\phantom{-}} E^{Y}_{s-m} 
\right)
\cr
&= \sum_{m=0}^{s-1} \frac{s-m}{s}\sum_{l=1}^{s-m} \frac{l}{s-m} 
   \left( L_{X_{l}} E^{X}_{s-m-l} E^{Y}_{m\phantom{-}} + 
          E^{X}_{m\phantom{-}} L_{Y_{l}} E^{Y}_{s-m-l} \right) 
\cr
&=\sum_{m=0}^{s-1}\sum_{l=1}^{s-m}\fraz{l}{s} 
   \biggl( L_{X_l}E^{X}_{s-m-l}E^{Y}_{m\phantom{-}} 
    +\sum_{k=0}^{m}L_{E^{X}_kY_l}E^{X}_{m-k}E^{Y}_{s-m-l}\biggr)
\cr
&= \sum_{l=1}^{s} \sum_{m=0}^{s-l} {l\over s} \left(
   L_{X_{l}} E^{X}_{s-m-l} E^{Y}_{m\phantom{-}} 
   +L_{Y_{l}} E^{X}_{m\phantom{-}} E^{Y}_{s-m-l} \right) 
\cr
&& +\sum_{l=1}^{s-1} \sum_{m=1}^{s-l} \sum_{k=1}^{m} {l\over s}
L_{E^{X}_k Y_{l}} E^{X}_{m-k} E^{Y}_{s-m-l}\ .\cr
}
$$
The identity~\frmref{trslie.5b} is used in order to obtain the third
line.  The first double sum in the latter expression is further
elaborated as
$$
\displaylines{
\qquad
\sum_{l=1}^{s} \frac{l}{s}\biggl(
 \lie{X_l} \sum_{m=0}^{s-l} E^{X}_{s-l-m} E^{Y}_m +
 \lie{Y_l} \sum_{m=0}^{s-l} E^{X}_{m} E^{Y}_{s-l-m}
  \biggr)
\hfill\cr\hfill
=\sum_{l=1}^{s} \frac{l}{s}\lie{X_l+Y_l} 
 \sum_{m=0}^{s-l} E^{X}_{m} E^{Y}_{s-l-m}
= \sum_{l=1}^{s} \frac{l}{s}\lie{X_l+Y_l} E^{Z}_{s-l}\ .
\qquad}
$$
The latter expression coincides with the first sum in the
r.h.s.~of~\frmref{trslie.8a}.  The triple sum is further elaborated as
$$
\displaylines{
\qquad
\sum_{l=1}^{s-1} \sum_{k=1}^{s-l} {l\over s} L_{E^{X}_k Y_{l}}
\sum_{h=0}^{s-l-k} E^{X}_h E^{Y}_{s-l-k-h}
=\sum_{l=1}^{s-1} \sum_{k=1}^{s-l} {l\over s} L_{E^{X}_k Y_{l}}
E^{Z}_{s-l-k} 
\hfill\cr\hfill
=\sum_{l=1}^{s-1} \sum_{m=l+1}^{s} {l\over s} L_{E^{X}_{m-l}Y_{l}}
E^{Z}_{s-m}
=\sum_{m=2}^{s} \sum_{l=1}^{m-1} {l\over s} L_{E^{X}_{m-l} Y_{l}}
E^{Z}_{s-m} \ .
\qquad\cr
}
$$
Here the induction hypothesis is used in the first step.  The last
expression coincides with the second double sum in the
r.h.s.~of~\frmref{trslie.8a}.  Thus, the right member
of~\frmref{trslie.8b} equals the last member of~\frmref{trslie.8a},
and this concludes the proof.

\references

\bye